\documentclass[12pt]{article}\usepackage{amssymb}

\newcommand{\text}{\textrm}

\newcommand{\AT}{{@}}
\newcommand{\cont}{{\frak c}}

\title{A model with no magic set}

\author{
Krzysztof Ciesielski%
\thanks{The authors wish to thank Professors Maxim Burke,
Andrzej Ros{\l}anowski and Jerzy Wojciechowski
for reading preliminary versions of this paper
and helping in improving its final version.
\endgraf
The work of the first author was partially supported by
NSF Cooperative
Reasearch Grant INT-9600548, with its Polish part being
financed by Polish Academy of Science PAN, and 1996/97 West Virginia
University Senate Research Grant. }
\\
{\footnotesize Department of Mathematics,}
{\footnotesize West Virginia University,}\\
{\footnotesize Morgantown, WV 26506-6310, USA}\\
{\footnotesize KCies\AT wvnvms.wvnet.edu}
\and
Saharon Shelah%
\thanks{This work was supported in part by a grant from
NSF grant DMS 97-04477.
Publication 653.
}\\
{\footnotesize Institute of Mathematics,}
{\footnotesize the Hebrew University of Jerusalem}\\
{\footnotesize 91904 Jerusalem, Israel}\\
{\footnotesize and}\\
{\footnotesize Department of Mathematics,}
{\footnotesize Rudgers University}\\
{\footnotesize New Brunswick, NJ 08854, USA}}
\date{}

\newcommand{\forces}{\mathrel{\|}\joinrel\mathrel{-}}
\newcommand{\suc}{{\rm succ}}
\newcommand{\ccc}{{\rm con}}

\newcommand{\nor}{{\rm norm}}
\newcommand{\norsup}{{\rm\underline{norm}}}
\newcommand{\sq}{\subseteq}
\newcommand{\real}{{\mathbb R}}

\newcommand{\rational}{{\mathbb Q}}

\newcommand{\la}{{\langle}}
\newcommand{\ra}{{\rangle}}
\newcommand{\restr}{{\hbox{$\,|\grave{}\,$}}}
\newcommand{\srestr}{{\hbox{${\scriptstyle\,|\grave{}\,}$}}}
\newcommand{\A}{{\cal A}}

\newcommand{\T}{{\cal T}}

 \newcommand{\Equi}{\mbox{$\Leftrightarrow$}}

\newcommand{\continuum}{{\cont}}

\def\poset{{\mathbb P}}

\def\dom{{\rm dom}}

\def\qed{\hfill\vrule height6pt width6pt depth1pt\medskip}

\newtheorem{theorem}{Theorem}[section]
\newtheorem{corollary}[theorem]{Corollary}
\newtheorem{proposition}[theorem]{Proposition}
\newtheorem{lemma}[theorem]{Lemma}
\newtheorem{problem}[theorem]{Problem}
\newtheorem{example}[theorem]{Example}
\newtheorem{definition}[theorem]{Definition}
\newtheorem{remark}[theorem]{Remark}
\newtheorem{Fact}[theorem]{Fact}

\newcommand{\thm}[2]{\begin{theorem}\label{#1}{\sl #2}\end{theorem}}
\newcommand{\cor}[2]{\begin{corollary}\label{#1}{\sl #2}\end{corollary}}
\newcommand{\prop}[2]{\begin{proposition}\label{#1}{\sl #2}\end{proposition}}
\newcommand{\lem}[2]{\begin{lemma}\label{#1}{\sl #2}\end{lemma}}

\newcommand{\rem}[2]{\begin{remark}\label{#1}{\rm #2}\end{remark}}

\begin{document}

\maketitle

\begin{abstract}
We will prove that there exists a model of ZFC+``$\continuum=\omega_2$''
in which every $M\sq\real$ of cardinality
less than continuum $\continuum$ is meager, and such that
for every $X\sq\real$ of cardinality $\continuum$
there exists a continuous function $f\colon\real\to\real$
with $f[X]=[0,1]$.

In particular in this model there is
no magic set, i.e., a set $M\sq\real$ such that
the equation $f[M]=g[M]$ implies $f=g$
for every continuous nowhere constant functions
$f,g\colon\real\to\real$.
\end{abstract}

\newpage

\section{Introduction}

The main goal of this paper is to prove the following theorem.

\thm{th:main}{There exists a model of ZFC in which $\continuum=\omega_2$,
\begin{description}
\item[($\star$)] for every $X\sq\real$ of cardinality $\continuum$
there exists a continuous function $f\colon\real\to\real$
such that $f[X]=[0,1]$, and

\item[($\star\star$)] every $M\sq\real$ of cardinality
less than $\continuum$ is meager.
\end{description}
}

Note that ($\star$) of Theorem~\ref{th:main} is known to hold
in the iterated perfect set model. (See A.~W.~Miller~\cite{Mi}.)
This result was also generalized by P.~Corazza~\cite{C} by finding
another model leading to the following theorem.

\thm{th:Corazza}{{\rm (Corazza)}
It is consistent with ZFC that ($\star$) holds and
\begin{description}
\item[($\star\star^\prime$)] every $M\sq\real$ of cardinality
less than $\continuum$ is of strong (so Lebesgue) measure zero.
\end{description}
}
Note that the condition ($\star\star$) is false in
the iterated perfect set model and in Corazza model.
(See~\cite{BC}.)

Corazza noticed also that Theorem~\ref{th:Corazza}
implies the following corollary (since there exists
a universal measure zero set of
cardinality ${\rm non}({\cal L})$, where
${\rm non}({\cal L})$ is the smallest cardinality
of a nonmeasurable set).
\cor{cor:Corazza}{{\rm (Corazza~\cite[Thm 0.3]{C})}
It is consistent with ZFC that ($\star$) holds and
there is a universal measure zero set of cardinality $\continuum$.
In particular in this model there are $2^\continuum$ many
universal measure zero sets of cardinality~$\continuum$.
}
He asked also whether the similar statement is true
with ``always first-category set'' replacing
``universal measure zero set.'' The positive answer
easily follows from Theorem~\ref{th:main}, since (in ZFC)
there exists
an always first-category set of
cardinality ${\rm non}({\cal M})$, where
${\rm non}({\cal M})$ is the smallest cardinality
of a nonmeager set.

\cor{cor:main0}{
It is consistent with ZFC that ($\star$) holds and
there is an always first-category set of cardinality $\continuum$.
In particular in this model there are $2^\continuum$ many
always first-category sets of cardinality~$\continuum$. \qed
}

Clearly Theorem~\ref{th:main}
can be viewed as dual to Theorem~\ref{th:Corazza}.
However, our original motivation for proving Theorem~\ref{th:main}
comes from another source.
In~\cite{BD} A.~Berarducci and D.~Dikranjan proved that
under the Continuum Hypothesis (abbreviated as
CH) there exists a set $M\sq\real$, called a {\em magic set},
such that for any two continuous nowhere constant functions
$f,g\colon\real\to\real$ if $f[M]\sq g[M]$ then $f=g$.
Different generalizations of a magic set were also studied
by M.~R.~Burke and K.~Ciesielski in~\cite{BC}. In particular
they examined the {\em sets of range uniqueness for the class $C(\real)$},
i.e., sets which definition is obtained from the definition
of a magic set by replacing the implication
``if $f[M]\sq g[M]$ then $f=g$'' with ``if $f[M]=g[M]$ then $f=g$.''
They proved \cite[Cor. 5.15 and Thm. 5.6(5)]{BC}
that if $M\sq\real$ is a
set of range uniqueness for $C(\real)$ then $M$ is not meager and
there is no continuous function $f\colon\real\to\real$ for which
$f[M]=[0,1]$.
This and Theorem~\ref{th:main} imply immediately
the following corollary, which solves the problems from
\cite{BD} and~\cite{BC}.

\cor{cor:main}{There exists a model of ZFC in which
there is no set of range uniqueness for $C(\real)$. In particular
there is no magic set in this model. }

Finally, it is worthwhile to mention that for the class of
nowhere constant differentiable function
the existence of a magic set is provable in ZFC, as noticed
by Burke and Ciesielski~\cite{BuCi2}.
In the same paper~\cite[cor.~2.4]{BuCi2} it has been noticed that
in the model constructed below there is also no
set or range uniquness
for $C(X)$ for any perfect Polish space $X$.

\section{Preliminaries}

Our terminology is standard and follows that from~\cite{BJ},
\cite{CiBook}, or \cite{Ku}.

A model satisfying Theorem~\ref{th:main}
will be obtained as a generic extension
of a model $V$ satisfying CH.
The forcing used to obtain such an extension
will be a countable support iteration $\poset_{\omega_2}$
of length $\omega_2$ of a forcing notion $\poset$ defined below.
Note that $\poset$, which is a finite level version of Laver forcing,%
\footnote{Note that Theorem~\ref{th:main} is false in Laver model,
since in this model there is a $\continuum$-Lusin set
(there is a scale) and such a set cannot be mapped continuously
onto $[0,1]$.}
is a version  of a tree-forcing
$\rational_1^{{\rm tree}}(K,\Sigma)$ from~\cite[sec.~2.3]{RS}
(for a 2-big finitary local tree-creating pair
$(K,\Sigma)$; it is also a relative of
the forcing notion defined in \cite[2.4.10]{RS})
and most of the results presented in this section
is a variation of general facts proved in this paper.
To define $\poset$, we need the following terminology.

A subset $T\sq\omega^{<\omega}$ is a {\em tree\/} if
$t\restr n\in T$ for every $t\in T$ and $n<\omega$.
For a tree $T\sq\omega^{<\omega}$ and $t\in T$ we will write
$\suc_T(t)$ for the set of all immediate successors of
$t$ in $T$, i.e.,
\[
\suc_T(t)=\{s\in T\colon t\sq s\ \&\ |s|=|t|+1\}.
\]
We will use the symbol $\T$ to denote the set of all
nonempty
trees $T\sq\omega^{<\omega}$ with no finite branches, i.e.,
\[
\T=\{T\sq\omega^{<\omega}\colon
\text{ $T\neq\emptyset$ is a tree\ \ \&\ \
     $\suc_T(t)\neq\emptyset$ for every $t\in T$}\}.
\]
For $T\in\T$ we will write $\lim T$ to denote
the set of all branches of $T$, i.e.,
\[
\lim T=\{s\in\omega^\omega\colon s\restr n\in T\text{ for every } n<\omega\}.
\]
Also if $t\in T\in\T$ then we define
\[
T^t=\{s\in T\colon s\sq t\ \text{ or }\ t\sq s\}.
\]

Now define inductively the following
``very fast increasing''
sequences
$\la b_i,n_i<\omega\colon i<\omega\ra$
by putting $n_{-1}=1$,
and for $i<\omega$
\[
b_i=(i+2)^{(n_{i-1}!)^i}
\ \ \ \text{ and }\ \ \
n_i=(b_i)^{(b_i)^i}.
\]
In particular $b_0=2$, $n_0=2$, $b_1=9$, $n_1=9^9$, $b_2=4^{[(9^9)!]^2}$,
etc. (For the purpose of our forcing any sequences that
grows at least ``as fast'' would suffice.)
Also let
\[
T^\star=\bigcup_{k<\omega}\prod_{i<k}n_i
=\{s\restr k\colon k<\omega\ \&\ s\in\prod_{i<\omega}n_i\}
\]
and
\[
\T^\star=\{T\in\T\colon T\sq T^\star\}.
\]
Forcing $\poset$ is defined as a family of all
trees $T\in\T^\star$ that have ``a lot of branching.''
To define this last term more precisely
we need the following definition
for every $i<\omega$, $T\in\T$ and $t\in T\cap \omega^i$:
\[
\nor_T(t)=
\log_{b_i}\log_{b_i}|\suc_T(t)|\in[-\infty,\infty).
\]
Note that $\nor_{T^\star}(t)=i$ for every
$t\in T^\star\cap \omega^i$.
Now for $T\in \T^\star$ and $k<\omega$
let
\[
\norsup_T(k)=\inf\{\nor_T(t)\colon t\in T\ \&\ |t|\geq k\}
\]
and define
\[
\poset=\left\{T\in\T^\star\colon
\lim_{k\to\infty}\norsup_T(k)=\infty
\right\}.
\]
The order relation on $\poset$ is standard. That is, $T_0\in\poset$
is stronger than $T_1\in\poset$, what we denote by $T_0\geq T_1$,
provided $T_0\sq T_1$.
Note also that $\norsup_T(k)\leq\norsup_{T^\star}(k)=k$
for every $k<\omega$.

In what follows for $t\in T\in\poset$
we will also use the following notation
\[
\norsup_T(t)=\norsup_{T^t}(|t|)
=\inf\{\nor_T(s)\colon s\in T\ \&\ t\sq s\}.
\]
It is easy to see that
\[
\norsup_T(k)=\min\{\norsup_T(t)\colon t\in T\cap\omega^k\}.
\]

For $n<\omega$ define a partial order $\leq_n$ on $\poset$
by putting $T_0\geq_n T$ if
\[
T_0\geq T\ \ \ \ \&\ \ \ \
T_0\restr\omega^k=T\restr\omega^k
\ \ \ \ \&\ \ \ \ \norsup_{T_0}(k)\geq n,
\]
where $k=\min\{j<\omega\colon\norsup_T(j)\geq n\}$.

Note that the sequence $\{\leq_n\colon n<\omega\}$
witnesses forcing $\poset$ to satisfy the axiom A.
(In particular $\poset$ is proper.)
That is (see \cite[7.1.1]{BJ} or \cite[2.3.7]{RS})
\begin{description}
\item{(i)} $T_0\geq_{n+1} T_1$ implies $T_0\geq_n T_1$
for every $n<\omega$ and $T_0,T_1\in\poset$;
\item{(ii)} if $\{T_n\colon n<\omega\}\sq\poset$ is such that
$T_{n+1}\geq_n T_n$ for every $n<\omega$
then there exists $T\in\poset$ extending each $T_n$, namely
$T=\bigcap_{n<\omega} T_n\in\poset$;
(such $T$ is often called a fusion of a sequence
$\la T_n\colon n<\omega\ra$;) and,

\item{(iii)} if $\A\sq\poset$ is an antichain,
then for every $T\in\poset$ and $n<\omega$ there exists
$T_0\in\poset$ such that $T_0\geq_n T$
and the set $\{S\in\A\colon S \text{ is compatible with } T\}$
is at most countable.
\end{description}
In fact, in case of the forcing $\poset$
the set $\{S\in\A\colon S \text{ is compatible with } T\}$
from (iii) is finite.
%
%
Since this fact will be heavily used in Section~\ref{sec3}
we will include here its proof. (See Corollary~\ref{corAxA}.)
However, this fact will not be used in the next three
sections so it can be skipped in the first reading.

The following definition is
a modification of the similar one for the Laver forcing.
(See~\cite[p. 353]{BJ}.)

\pagebreak

Let $D\sq\poset$ be dense below $p\in\poset$ and $n<\omega$.
For $t\in p$ with $\norsup_p(t)\geq n$
we define the ordinal number $r^n_D(t)<\omega$ as follows:\label{page:Rank}
\begin{description}
\item[(1)] $r^n_D(t)=0$ if there exists $p^\prime\in D$ extending $p^t$
           such that $\norsup_{p^\prime}(t)\geq n-1$;
\item[(2)] if $r^n_D(t)\neq 0$ and $t\in\omega^i$ then
           \[
           r^n_D(t)=\min\left\{\alpha\colon
           \left(\exists U\in[\suc_p(t)]^{\geq (b_i)^{(b_i)^{n-1}}}\right)
           \bigl(\forall s\in U\bigr)\left(r^n_D(s)<\alpha\right)\right\}.
           \]
\end{description}

\lem{lem0xx}{Let $D\sq\poset$ be dense below $p\in\poset$ and
$n<\omega$. Then $r^n_D(t)$ is well defined for every
$t\in p$ with $\norsup_p(t)\geq n$.}

Proof. By way of contradiction assume that there exists
$t\in p$ with $\norsup_p(t)\geq n$ for which $r^n_D(t)$ is undefined.
Then $n>1$ (since otherwise we would have $r_D^n(t)=0$)
and for any such $t$ belonging to $\omega^i$ the set
$$
U=\{s\in\suc_p(t)\colon r^n_D(s)\text{ is defined}\}
$$
has cardinality less than $(b_i)^{(b_i)^{n-1}}$.
So 
\begin{equation}\label{con1}
\!\!\!\!\!\!
\!\!\!
|\{s\in\suc_p(t)\colon r^n_D(s)\text{ is undefined}\}|=
|\suc_p(t)\setminus U|\geq |\suc_p(t)|/2
\end{equation}
since $|\suc_p(t)|/2=(b_i)^{(b_i)^{\nor_p(t)}}/2
\geq (b_i)^{(b_i)^n}/2\geq (b_i)^{(b_i)^{n-1}}>|U|$.
Construct a tree $p_0\in\T^\star$ such that  $p_0\sq p^t$,
\begin{equation}\label{con2}
|\suc_{p_0}(s)|\geq|\suc_p(s)|/2
\end{equation}
and $r^n_D(s)$ is undefined for every $s\in p_0$ with $t\sq s$.
The construction can be easily done
by induction on the levels of a tree,
using (\ref{con1}) to make an inductive step.
But
(\ref{con2}) implies that
for every $i<\omega$ and $s\in p_0\cap\omega^i$ with $t\sq s$
$$
\nor_{p_0}(s)=\log_{b_i}\log_{b_i}|\suc_{p_0}(s)|\geq
\log_{b_i}\log_{b_i}|\suc_p(s)|/2
\geq\nor_p(s)-1.
$$
So $p_0\in\poset$. Take $p^\prime\in D$
with $p^\prime\geq p_0$. We can find
$t_1\in p^\prime$
such that $\norsup_{p}(t_1)\geq\norsup_{p^\prime}(t_1)\geq n-1$.
Then $r^n_D(t_1)=0$, contradicting the fact that
$r^n_D(s)$ is undefined for every $s\in p_0\supseteq p^\prime$.
\qed

\lem{lem00xx}{
Let $D\sq\poset$ be dense below $p\in\poset$ and
$n<\omega$. Then for every
$t\in p$ with $\norsup_p(t)\geq n$ there exist $p_t\geq_{n-1} p^t$
and a finite set $A_t\sq p_t$ such that
$p_t=\bigcup_{s\in A_t}(p_t)^s$ and $(p_t)^s\in D$ for every $s\in A_t$.
}

Proof. The proof is by induction on $r^n_D(t)$.

If $r^n_D(t)=0$ then $p_t=p^\prime\in D$
will satisfy the lemma with $A_t=\{t\}$.

If $r^n_D(t)=\alpha>0$ choose $U\in[\suc_p(t)]^{\geq (b_i)^{(b_i)^{n-1}}}$
from the definition of $r^n_D(t)$.
By the inductive assumption for every $u\in U$ there exist
$q_u\geq_{n-1}p^u$ and a finite set $A_u\sq q_u$ such that
$q_u=\bigcup_{s\in A_u}(q_u)^s$ and $(q_u)^s\in D$ for every
$s\in A_u$.
Then $p_t=\bigcup_{u\in U}q_u$ and $A_t=\bigcup_{u\in U}A_u$
satisfy the lemma.
\qed

The next corollary can be also found, in general form,
in~\cite[2.3.7, 3.1.1]{RS}.

\cor{corAxA}{Let $\A\sq\poset$ be an antichain.
Then for every $p\in\poset$ and $n<\omega$ there exists
$q\in\poset$ such that $q\geq_n p$
and the set
\[
\A_0=\{r\in\A\colon r \text{ is compatible with } q\}
\]
is finite.
}

Proof. Extending $\A$, if necessary, we can assume that $\A$ is a
maximal antichain. Thus
$D=\{q\in\poset\colon (\exists p\in\A)(q\geq p)\}$
is dense in $\poset$.

Let $i<\omega$ be such that $\norsup_p(i)\geq n+1$.
By Lemma~\ref{lem00xx} for every $t\in p\cap\omega^i$
there exists $p_t\geq_n p^t$
and a finite set $A_t\sq p_t$ such that
$p_t=\bigcup_{s\in A_t}(p_t)^s$ and $(p_t)^s\in D$ for every $s\in A_t$.
Put
$q=\bigcup_{t\in p\cap\omega^i}p_t$.
Then it satisfies the corollary. \qed

\section{Proof of the theorem}

For $\alpha\leq\omega_2$ let $\poset_\alpha$ be a countable
support iteration of forcing $\poset$ defined in the previous section.
Thus $\poset_\alpha$ is obtained from a sequence
$\la\la\poset_\beta,\dot\rational_\beta\ra\colon\beta<\alpha\ra$,
where each $\poset_\beta$ forces that $\dot\rational_\beta$ is a
$\poset_\beta$-name for forcing $\poset$. Also we will
consider elements of $\poset_\alpha$ as functions
$p$ which domains are countable subset of $\alpha$.
In particular if $p\in\poset_\alpha$ and $0\in\dom(p)$
then $p(0)$ is an element of $\poset$ as defined in $V$.

Now let $V$ be a model of ZFC+CH and let $G$ be a $V$-generic
filter in $\poset_{\omega_2}$. We will show that the conclusion of
Theorem~\ref{th:main} holds in $V[G]$.

In what follows for $\alpha\leq\omega_2$
we will use the symbol $G_\alpha$ to denote
$G\cap\poset_\alpha$. In particular each
$G_\alpha$ is a $V$-generic filter in $\poset_\alpha$ and
$V[G_\alpha]\sq V[G_{\omega_2}]=V[G]$.

Since CH holds in $V$, forcing $\poset_{\omega_2}$ is $\omega_2$-cc in $V$.
Thus since $\poset$ satisfies the axiom A,
we conclude that $\poset_{\omega_2}$ preserves cardinal numbers
and indeed $\continuum=\omega_2$ holds in $V[G]$.

To prove that $(\star\star)$ holds in $V[G]$
consider $\prod_{i<\omega}n_i=\lim T^\star$
with the product topology. Since $\prod_{i<\omega}n_i$
is homeomorphic to the Cantor set $2^\omega$ it is enough to show that
every subset $S$ of $\prod_{i<\omega}n_i$ of
cardinality less than $\continuum^{V[G]}=\omega_2$
is meager in $\prod_{i<\omega}n_i$. But every $x\in S$ belongs
already to some intermediate model $V[G_\alpha]$
with $\alpha<\omega_2$, since $\poset$ satisfies the axiom A
(so is proper), and the iteration is with countable support.
In particular there
exists an
$\alpha<\omega_2$ such that
$S\sq V[G_\alpha]$. So it is enough to prove that
$\left(\prod_{i<\omega}n_i\right)\cap V[G_\alpha]$
is meager in $\prod_{i<\omega}n_i$.

Since $V[G_{\alpha+1}]$ is obtained from $V[G_\alpha]$
as a generic extension via forcing $\poset$ (in $V[G_\alpha]$)
our claim concerning $(\star\star)$ in $V[G]$
follows immediately from the following lemma.
(See also \cite[3.2.8]{RS}.)

\lem{lem2star}{Let $V$ be a model of ZFC+CH and
$H$ be a $V$-generic filter in $\poset$.
Then in $V[H]$ the set
$\left(\prod_{i<\omega}n_i\right)\cap V$
is a meager subset of $\prod_{i<\omega}n_i$.}

Proof. Let $r\in\prod_{i<\omega}n_i$ be such that
$\{r\}=\bigcap\{\lim T\colon T\in H\}$
and put $M=\bigcup_{j<\omega}M_j$ where
\[
M_j=\left\{s\in\prod_{i<\omega}n_i\colon s(k)\neq r(k)
\text{ for every } j\leq k<\omega\right\}.
\]
Since clearly every $M_j$ is closed nowhere dense it is
enough to show that $\left(\prod_{i<\omega}n_i\right)\cap V\sq M$.
For this pick $s\in \left(\prod_{i<\omega}n_i\right)\cap V$
and consider a subset $D=\bigcup_{j<\omega}D_j\in V$ of $\poset$,
where
\[
D_j=\{p\in\poset\colon
(\forall t\in p)(\forall k\in\dom(t)\setminus j)(s(k)\neq t(k))\}.
\]
It is enough to prove that $D$ is dense in $\poset$,
since $H\cap D_j\neq \emptyset$ implies that $s\in M_j$.

So let $p_0\in\poset$ and let $j<\omega$ be such that
$p_0\cap\omega^{j-1}\sq p_0(1)$ and define
\[
p=\{t\in p_0\colon (\forall k\in\dom(t)\setminus j)(s(k)\neq t(k))\}.
\]
Clearly $p$ is a tree.
It is enough to show that $p\in\poset$, since then $p\in D_j$ extends $p_0$.
But if $t\in p_0\cap\omega^k$ for some $k\geq j$
and $t_0\in p$ is an immediate predecessor of $t$ then
\[
|\suc_p(t_0)|\geq|\suc_{p_0}(t_0)|-1
=(b_{k-1})^{(b_{k-1})^{\nor_{p_0}(t_0)}}-1>0
\]
so $\suc_p(t_0)$ is nonempty and for every $s\in\lim p$
\[
\lim_{i\to\infty}\norsup_p(s\restr i)\geq
\lim_{i\to\infty}(\norsup_{p_0}(s\restr i)-1)
=\infty.
\]
This finishes the proof of Lemma~\ref{lem2star}. \qed

To show that $(\star)$ holds in $V[G]$ we will use the following two
propositions.
The first of them is an easy modification of the
Factor Theorem from \cite[Thm~1.5.10]{BJ}.
For the case of Sacks forcing
this has been proved in~\cite[Thm~2.5]{BL}.

\prop{prop1}{Let $\beta<\alpha\leq\omega_2$ and $\gamma$
be such that $\beta+\gamma=\alpha$.
If $\poset^\star_\gamma$ is a $\poset_\beta$-name
for the iteration $\poset_\gamma$ of $\poset$
(as constructed in $V^{\poset_\beta}$) then
forcings $\poset_\alpha$ and $\poset_\beta\star\poset^\star_\gamma$
are equivalent.\qed}


The analog of the next proposition for the iteration of Sacks forcing
can be found in an implicit form in~\cite{Mi}.

\prop{prop2}{Suppose that
$p\forces``\tau\in 2^\omega\setminus V$''
for some $p\in\poset_{\omega_2}$. Then (in $V$) there exists
a continuous function $f\colon 2^\omega\to 2^\omega$ with the property that
\begin{itemize}
\item for every $r\in 2^\omega$ there exists $q_r\geq p$ such that
\[
q_r\forces f(\tau)=r.
\]
\end{itemize}
}

The proof of
Proposition~\ref{prop2} will be postponed to the next section.
The proof of $(\star)$
based on Proposition~\ref{prop2} and
presented below is an elaboration of the proof
from~\cite{Mi} that $(\star)$ holds in the iterated Sacks model.

First note (compare \cite{C})
that to prove $(\star)$ it is enough to show that
\begin{description}
\item[$(\circ)$]
for every $X\sq 2^\omega$ of cardinality $\continuum$
there exists a continuous function
$f\colon 2^\omega\to 2^\omega$ such that
$f[X]= 2^\omega$.
\end{description}

Indeed if $X\sq\real$ has cardinality $\continuum$
and there is no zero-dimensional perfect set $P\subset\real$
such that $|X\cap P|=\continuum$ then $X$ is a
$\continuum$-Lusin subset of $\real$.
Then there is $\continuum$-Lusin subset of $ 2^\omega$
as well, and such a set would contradict $(\circ)$
since it cannot be mapped continuously onto $[0,1]$
(so onto $2^\omega$ as well).
(See e.g. \cite[Sec.~2]{Mi}.)

So there are $a,b\in\real$ and a
zero-dimensional perfect set $P\subset[a,b]$
with $|X\cap P|=\continuum$. But $P$ and
$ 2^\omega$ are homeomorphic. Therefore, by $(\circ)$,
there exists a continuous $f\colon P\to P\subset[a,b]$ such that
$f[X\cap P]=P$. Then a continuous extension
$F\colon\real\to[a,b]$ of $f$, which exists by Tietze Extension theorem,
has a property that $F[X]\supseteq P$. Now if $g\colon\real\to[0,1]$
is continuous and
such that $g[P]=[0,1]$ then $f=g\circ F$ satisfies $(\star)$.

We will prove $(\circ)$ in $V[G]$ by contraposition.
So let $X\sq 2^\omega$ be such that
$f[X]\neq 2^\omega$ for every continuous
$f\colon 2^\omega\to 2^\omega$.
Thus for any such $f$ there exists an $F_0(f)\in 2^\omega$
such that $F_0(f)\notin f[X]$.
We will prove that this implies $|X|<\continuum$
by showing that $X\sq V[G_\alpha]$ for some $\alpha<\omega_2$.
This is enough, since $V[G_\alpha]$ satisfies CH.

Now let $D=2^{<\omega}\in V$. Since $D$ is dense in $2^\omega$
any continuous $f\colon 2^\omega\to 2^\omega$ is uniquely
determined by $f\restr D$. Let $F\colon (2^\omega)^D\to 2^\omega$,
$F\in V[G]$, be such that $F(f\restr D)=F_0(f)$
for every continuous $f\colon 2^\omega\to 2^\omega$. Thus
\[
F(f\restr D)\notin f[X]
\]
for every continuous $f\colon 2^\omega\to 2^\omega$.
We claim that there exists an $\alpha<\omega_2$
of cofinality $\omega_1$ such that
\begin{equation}\label{thmeq1}
F\restr\left((2^\omega)^D\cap V[G_\alpha]\right)\in  V[G_\alpha].
\end{equation}

To show (\ref{thmeq1}) first
recall that for every real number $r$
and every $\alpha\leq\omega_2$ of uncountable cofinality
if $r\in V[G_\alpha]$ then $r\in V[G_\beta]$
for some $\beta<\alpha$. This is a general property of a countable
support iteration of forcings satisfying the axiom A (and, more generally,
proper forcings). In particular
\begin{equation}\label{eq2}
(2^\omega)^D\cap V[G_\alpha]=
\bigcup_{\beta<\alpha}\left((2^\omega)^D\cap V[G_\beta]\right)
\end{equation}
for every $\alpha<\omega_2$ of cofinality $\omega_1$.

Now let $\la f_\alpha\colon\alpha<\omega_2\ra\in V[G]$
be a one-to-one enumeration of $(2^\omega)^D$, and put
$y_\alpha=F(f_\alpha)$. Then there exists a sequence
$S=\la\la\varphi_\alpha,\eta_\alpha\ra\colon \alpha<\omega_2\ra\in V$
such that $\varphi_\alpha$ and $\eta_\alpha$ are the
$\poset_{\omega_2}$-names for $f_\alpha$ and $y_\alpha$, respectively.
Moreover, since $\poset_{\omega_2}$ is $\omega_2$-cc in $V$,
we can assume that for every $\alpha<\omega_2$
there is a $\delta(\alpha)<\omega_2$ such that
$\varphi_\alpha$ and $\eta_\alpha$ are the
$\poset_{\delta(\alpha)}$-names.
Also if we choose $\delta(\alpha)$ as the
smallest number with this property, then
function $\delta$ belongs to $V$, since it is
definable from $S\in V$.

Note also that for every $\beta<\omega_2$ there is an
$h_0(\beta)<\omega_2$ with the property that
for every $f\in (2^\omega)^D\cap V[G_\beta]$
there is $\gamma<h_0(\beta)$ such that $\varphi_\gamma$
is a name for $f$ (with respect to $G$).
Once again using the fact that $\poset_{\omega_2}$ is $\omega_2$-cc
in $V$ we can find in $V$ a function $h\colon\omega_2\to\omega_2$
bounding $h_0\in V[G]$, i.e.,
such that $h_0(\beta)\leq h(\beta)$ for every $\beta<\omega_2$.
Let
\[
C=\{\alpha<\omega_2\colon
(\forall\gamma<\alpha)(\delta(\gamma),h(\gamma)<\alpha)\}\in V.
\]
Then $C$ is closed and unbounded in $\omega_2$. Pick $\alpha\in C$
of cofinality $\omega_1$. We claim that $\alpha$ satisfies (\ref{thmeq1}).

To see it, note first that the definition of $\delta$
implies that every name in the sequence
$\la\la\varphi_\gamma,\eta_\gamma\ra\colon \gamma<\alpha\ra$
is a $\poset_\alpha$-name. So
\[
F\restr\{f_\gamma\colon\gamma<\alpha\}=
\{\la f_\gamma,y_\gamma\ra\colon \gamma<\alpha\}\in V[G_\alpha].
\]
Moreover clearly
$\{f_\gamma\colon\gamma<\alpha\}\sq (2^\omega)^D\cap V[G_\alpha]$.
However, by (\ref{eq2}),
for every $f\in (2^\omega)^D\cap V[G_\alpha]$
there exists $\beta<\alpha$ such that
$f\in (2^\omega)^D\cap V[G_\beta]$.
Thus, by the definition of $h_0$ and $h$, there exists
$\gamma<h_0(\beta)\leq h(\beta)<\alpha$ such that
$f=f_\gamma$. So
$\{f_\gamma\colon\gamma<\alpha\}=(2^\omega)^D\cap V[G_\alpha]$
and (\ref{thmeq1}) has been proved.

Now take an $\alpha<\omega_2$ having property 
(\ref{thmeq1}).
For this $\alpha$ we will argue that $X\sq V[G_\alpha]$.
But, by Proposition~\ref{prop1},
$V[G]$ is a generic extension of $V[G_\alpha]$
via forcing $\poset_{\omega_2}$ as defined in
$V[G\cap\poset_\alpha]$.
Thus without loss of generality we can assume that $V[G_\alpha]=V$.
In particular
\[
F_1=F\restr\left((2^\omega)^D\cap V\right)\in  V.
\]

To see that $X\sq V$ take an arbitrary
$z\in 2^\omega\setminus V$, and
pick a $\poset_{\omega_2}$-name $\tau$ for $z$.  Let
$p_0\in G\subset\poset_{\omega_2}$ be such that
$p_0\forces``\tau\in 2^\omega\setminus V$''
and fix an arbitrary $p_1\geq p_0$.
Working in $V$ we will find a $p\in\poset_{\omega_2}$ stronger than $p_1$
and a continuous function $f\in V$ from $2^\omega$ to $2^\omega$
such that
\begin{equation}\label{eqThm}
p\forces f(\tau)=F_1(f\restr D).
\end{equation}

To see it notice that by Proposition~\ref{prop2} there exists
a continuous function $f\colon 2^\omega\to 2^\omega$ such that
for every $r\in 2^\omega$ (from $V$) there exists $q_r\geq p_1$ with
\[
q_r\forces f(\tau)=r.
\]
Take $r=F_1(f\restr D)\in V$. Then $p=q_r$ satisfies (\ref{eqThm}).

Now (\ref{eqThm}) implies that the set
\[
E=\{q\in\poset_{\omega_2}\colon
(\exists\text{ continuous }f\colon 2^\omega\to 2^\omega)
(q\forces``f(\tau)=F_1(f\restr D))\text{''}\}\in V
\]
is dense above $p_0\in G$.
Therefore, there exist $q\in G\cap E$ and a
continuous function $f\colon 2^\omega\to 2^\omega$ such that
$q\forces``f(\tau)=F_1(f\restr D)$.''
In particular $f(z)=F_1(f\restr D)=F(f\restr D)\notin f[X]$,
implying that $z\not\in X$.
Since it is true for every $z\in 2^\omega\setminus V$,
we conclude that $X\subset V$.

This finishes the proof of Theorem~\ref{th:main}
modulo the proof of Proposition~\ref{prop2}.
\qed

\section{Proof of Proposition~\ref{prop2} --- another reduction}

In this short section we will prove Proposition~\ref{prop2}
based on one more technical lemma. The proof of the lemma
will be postponed to the next section.

To state the lemma and prove the proposition
we need the following iteration version of the axiom A.
For $\alpha\leq\omega_2$, $F\in[\alpha]^{<\omega}$, and $n<\omega$ define
a partial order relation $\leq_{F,n}$ on $\poset_\alpha$
by
\[
q\geq_{F,n} p\ \ \Equi\ \ q\geq p\ \&\
(\forall\xi\in F)(q\restr\xi\forces q(\xi)\geq_n p(\xi)).
\]
Note that if $\xi\notin\dom(p)$ for some $\xi\in F$ then
it might be unclear what we mean by $p(\xi)$ in the above
definition. However, in such a case we will
identify $p$ with its extension, for which we put
$p(\xi)=\hat T^\star$, where $\hat T^\star$ is the standard
$\poset_\xi$-name for the weakest element $T^\star$ of $\poset$.
Recall also that if an increasing sequence $\la F_n\colon n<\omega\ra$
of finite subsets of $\alpha$ and
$\la p_n\in\poset_\alpha\colon n<\omega\ra$
are such that $p_{n+1}\geq_{F_n,n} p_n$ for every $n<\omega$
and $\bigcup_{n<\omega}\dom(p_n)=\bigcup_{n<\omega}F_n$ then
there exists $q\in\poset_\alpha$ extending each $p_n$.
(See e.g. \cite[7.1.3]{BJ}.)

\lem{lem:Split}{Let $\alpha<\omega_2$,
$p\in\poset_\alpha$ and $\tau$ be a $\poset_\alpha$-name such that
for every $\gamma<\alpha$
\[
p\forces \tau\in 2^\omega\cap V[G_\alpha]\setminus V[G_\gamma].
\]
Then there exists $q_0\in\poset_\alpha$ stronger than $p$
such that for
every $F\in[\alpha]^{<\omega}$, $n<\omega$, and
$q\in\poset_\alpha$ extending $q_0$
there exist (in V) an $m<\omega$, nonempty disjoint sets
$B_0,B_1\subset 2^m$ and $p_0,p_1\geq_{F,n} q$
such that
\[
p_j\forces\tau\restr m\in B_j
\]
for $j<2$.
}

Basically Lemma~\ref{lem:Split} is true since $p$ forces that
$\tau$ is a new real number. However, its proof
is quite technical and will be postponed
for the next section.

Next we will show
how Lemma~\ref{lem:Split} implies Proposition~\ref{prop2}.

\medskip

\noindent
{\bf Proof of Proposition~\ref{prop2}.} Let
$p\in\poset_{\omega_2}$ and $\tau$ be a $\poset_{\omega_2}$-name such that
$p\forces``\tau\in 2^\omega\setminus V$.''
Then, replacing $p$ with some stronger condition if necessary,
we can assume that there exists $\alpha<\omega_2$
such that for every $\gamma<\alpha$
\[
p\forces \tau\in V[G_\alpha]\setminus V[G_\gamma].
\]
In particular, since $p\forces``\tau\in V[G_\alpha]$,''
we can assume that $\tau$ is a $\poset_\alpha$-name.
We can also find $p^\prime\geq p\restr\alpha$ such that
\[
p^\prime\forces \tau\in 2^\omega\cap V[G_\alpha]\setminus V[G_\gamma].
\]
Thus it is enough to assume that $p\in\poset_\alpha$
and find $f\in V$ and
$q_r\in\poset_\alpha$ satisfying Proposition~\ref{prop2}.
(Otherwise, we can replace $q_r$'s with
$q_r\cup p\restr(\omega_2\setminus\alpha)$.)

For $m<\omega$ and $B\sq 2^m$ let
$[B]=\{x\in 2^\omega\colon x\restr m\in B\}$. Thus $[B]$ is a clopen
subset of $2^\omega$.
For every $s\in 2^{<\omega}$ we will define
$q_s\in\poset_\alpha$, $m_s<\omega$ and $B_s\sq 2^{m_s}$.
The construction will be done by induction on length $|s|$ of $s$.
Simultaneously we will construct an increasing sequence
$\la F_n\in[\alpha]^{<\omega}\colon n<\omega\ra$
such that the following conditions are satisfied
for every $s\in 2^{<\omega}$ and $n=|s|$:
\begin{description}
\item[(I0)]
$\bigcup\{\dom(q_t)\colon t\in 2^{<\omega}\}=\bigcup_{n<\omega}F_n$;
\item[(I1)] $q_{s0},q_{s1}\geq_{F_n,n} q_s$;
\item[(I2)] $B_{s0}\cap B_{s1}=\emptyset$, and
            $[B_{s0}]\cup [B_{s1}]\sq [B_s]$;
\item[(I3)] $q_{sk}\forces``\tau\restr m_s\in B_{sk}$'' for every $k<2$.
\end{description}

It is easy to fix an inductive schema of choice of
$F_n$'s which will force condition (I0) to be satisfied.
Thus we will assume that we are using such a schema
throughout the construction, without specifying its details.

Now let $q_0$ be as in Lemma~\ref{lem:Split}.
This will be our $q_\emptyset$.
Moreover if $q_s$ is already defined for some $s\in 2^{<\omega}$
then we choose $m_s$, $q_{s0}$, $q_{s1}$, $B_{s0}$, and $B_{s1}$
by using  Lemma~\ref{lem:Split} for $q=q_s\geq q_0$, $n=|s|$ and $F=F_n$.
This finishes the inductive construction.

Next for $n<\omega$ let $m_n=\max\{m_s\colon s\in 2^{\leq n}\}$
and for $s\in 2^n$ and $k<2$ put
$B_{sk}^\star=\{t\in 2^{m_n}\colon t\restr m_s\in B_{sk}\}$.
Then $p_{sk}\forces``\tau\restr m_n\in B_{sk}^\star$'' for every $k<2$.
Thus, replacing sets $B_{sk}$ with $B_{sk}^\star$ if necessary,
we can assume that $m_s=m_n$ for every $s\in 2^n$.

Note also that $\lim_{n\to\infty} m_n=\infty$. This follows easily
from (I3) and (I2).
Let
\[
P=\bigcap_{n<\omega}\bigcup_{s\in 2^n}[B_s].
\]
Then $P$ is perfect subset of $2^\omega$.
Define function $f_0\colon P\to 2^\omega$
by putting
\[
f_0(x)=r\ \ \text{ if and only if }\ \
x\in [B_{r\srestr n}] \text{ for every $n<\omega$}.
\]
It is easy to see that $f_0$ is continuous. Thus, by Tietze Extension
theorem, we can find a continuous extension
$f\colon 2^\omega\to 2^\omega$ of $f_0$.
We will show that $f$ satisfies the requirements
of Proposition~\ref{prop2}.

Indeed take $r\in 2^\omega$ and let
$q_n=q_{r\srestr n}$. Then, by (I1), $q_{n+1}\geq_{F_n,n} q_n$
for every $n<\omega$. Moreover, by (I0),
$\bigcup_{n<\omega}\dom(q_n)\sq \bigcup_{n<\omega}F_n$.
In addition, we can assume that the equation holds,
upon the identification described in the definition of
$\geq_{F,n}$. Thus there exists a $q_r\in\poset_\alpha$
extending each $q_n$.
But for every $n<\omega$
\[
q_{n+1}\forces\tau\restr m_n\in B_{r\srestr n+1}
\]
so that
\[
q_r\forces f_0( [\{\tau\restr m_n\}]\cap P )\in
f_0([B_{r\srestr n+1}]\cap P)\sq [\{r\restr n+1\}].
\]
Therefore, by the continuity of $f$,
\[
q_r\forces f(\tau)=r.
\]
This finishes the proof of Proposition~\ref{prop2}.


\section{Proof of Lemma~\ref{lem:Split}}\label{sec3}

We will start this section with the following
property that will be used several times in the sequel.

\lem{cor00}{Forcing $\poset$ has the property \text{\bf B}
from~\cite[p. 330]{BJ}. That is, for every $p\in\poset$,
a $\poset$-name $\mu$,
and $k<\omega$, if $p\forces``\mu\in\omega$''
then there exist $m<\omega$ and $p^\prime\geq_k p$ such that
$p^\prime\forces``\mu\leq m$.''}

Proof. This follows immediately from Corollary~\ref{corAxA}
applied to a maximal antichain
in the set $D=\{q\geq p\colon (\exists m<\omega)(q\forces\mu=m)\}$.
\qed


Recall also the following result concerning the property~\text{\bf B}.
(See~\cite[Lemma~7.2.11]{BJ}.)

\cor{Cor:lem11}{Let $\alpha\leq\omega_2$. If $p\in\poset_\alpha$,
$n\in\omega$, $F\in[\omega_2]^{<\omega}$
and $p\forces``\mu\in\omega$''
then there exist $m<\omega$ and $p^\prime\geq_{F,n} p$ such that
$p^\prime\forces``\mu<m$.'' }

The difficulty of the proof of  Lemma~\ref{lem:Split}
comes mainly from the
fact that  we have to find ``real'' sets $B_0$ and $B_1$
using for this only
$\poset_\alpha$-name $\tau$, and $p\in\poset_\alpha$,
which is also formed mainly from different names.
For this we will have to describe how to
recover ``real pieces of information''
form $\tau$ and $p$.
We will start this with the following lemma.

\lem{lem:ReadCon1}{Let $p\in\poset$ and $\tau$ be a $\poset$-name
such that $p\forces``\tau\in 2^\omega$.''
Then for every $n,m<\omega$ there exist $q\geq_n p$,
$i<\omega$, and a family
$\{x_s\in 2^m\colon s\in q\cap\omega^i\}$
such that for any $s\in q\cap\omega^i$
\[
q^s\forces \tau\restr m=x_s.
\]
}

Proof. Let
$D=\{p\in\poset\colon(\exists x\in 2^m)(p\forces\tau\restr m=x)\}$
and let $j<\omega$ be such that $\norsup_p(j)\geq n+1$.
By Lemma~\ref{lem00xx} for every $t\in p\cap\omega^j$
there exist $p_t\geq_n p^t$
and a finite set $A_t\sq p_t$ such that
$p_t=\bigcup_{s\in A_t}(p_t)^s$ and $(p_t)^s\in D$ for every $s\in A_t$.
Put
$q=\bigcup_{t\in p\cap\omega^j}p_t$
and let $i<\omega$ be such that
$\bigcup\{A_t\colon t\in p\cap\omega^j\}\sq\omega^{\leq i}$.
Then $q$ and $i$ satisfy the requirements. \qed

Let $p\in\poset$ and $\tau$ be a $\poset$-name
such that $p\forces``\tau\in 2^\omega$.''
We will say that {\em $p$ reads $\tau$ continuously\/}
if for every $m<\omega$ there exist
$i_m<\omega$ and a family
$\{x_s\in 2^m\colon s\in p\cap\omega^{i_m+1}\}$
such that for any $s\in p\cap\omega^{i_m+1}$
\[
p^s\forces \tau\restr m=x_s.
\]

\lem{lem:ReadCon2}{Let $p\in\poset$ and $\tau$ be a $\poset$-name
such that $p\forces``\tau\in 2^\omega$.''
Then for every $n<\omega$ there exists $q\geq_n p$
such that $q$ reads $\tau$ continuously.
}

Proof. By Lemma~\ref{lem:ReadCon1} we can define inductively
a sequence $\la q_m\colon m<\omega\ra$
such that $q_0=p$ and for every $m<\omega$
\begin{itemize}
\item $q_{m+1}\geq_{n+m} q_m$; and,
\item there exist $i_m<\omega$ and a family
$\{x_s\in 2^m\colon s\in q\cap\omega^{i_m+1}\}$
such that for any $s\in q\cap\omega^{i_m+1}$
\[
q^s\forces \tau\restr m=x_s.
\]
\end{itemize}
Then the fusion $q=\bigcap_{m<\omega} q_m$ of all $q_m$'s
has the desired properties. \qed

The next lemma is an important step in
our proof of Lemma~\ref{lem:Split}.
It also implies it quite easily
for $\poset_\alpha=\poset$. (See Corollary~\ref{lem:P0cor}.)

\lem{lem:P0}{Let $p\in\poset$ and $\tau$ be a $\poset$-name
such that
\[
p\forces\tau\in 2^\omega\setminus V
\]
and $p$ reads $\tau$ continuously with the sequence
$\la i_m\colon m<\omega\ra$ witnessing it.
Then for every $n,k<\omega$ with $\norsup_p(k)\geq n+1\geq 2$
there exist an arbitrarily large number $m<\omega$
and $q\geq_n p$ which can be represented as
\begin{equation}\label{eqZZ}
q=\bigcup_{t\in A}p_t,
\end{equation}
where
$A\sq p\cap\left(\omega^{\leq i_m}\setminus\omega^{<k}\right)$ and
the elements of $A$ are pairwise incompatible (as functions).
Moreover
for every $t\in A$ we have $p_t\geq p^t$,
and there exists a \underline{one-to-one} mapping
$\suc_q(t)\ni s\longmapsto x_s\in 2^m$
such that
\[
q^s\forces\tau\restr m=x_s
\]
for every $s\in\suc_q(t)$.
}

Proof. Fix $p$, $\tau$, $n$, and $k$ as in the lemma.
For every $u\in p\cap\omega^k$ and $m<\omega$ with $i=i_m>k$
consider the following {\em trimming procedure}.

For every $s\in p^u\cap\omega^{i+1}$ let
$x_s\in 2^m$ be such that
\[
p^s\forces\tau\restr m=x_s,
\]
put $q_{i+1}=p^u$, and assign to every  $s\in p^u\cap\omega^{i+1}$
a tag ``constant $x_s$.''
By induction we define a sequence
\[
q_{i+1}\leq q_i\leq q_{i-1}\leq q_{i-2}\leq \ldots\leq q_k
\]
of elements of $\poset$ such that for $k\leq j\leq i$
every $t\in q_j\cap\omega^j$
has a tag of either ``one-to-one'' or a ``constant $x_t$''
with $x_t\in 2^m$.

If for some $j\geq k$ the tree $q_{j+1}$
is already defined then for every $t\in q_{j+1}\cap\omega^j$
choose $U_t\sq\suc_{q_{j+1}}(t)$ of cardinality
$\geq .5\ |\suc_{q_{j+1}}(t)|\geq |\suc_{q_{j+1}}(t)|^{1/2}$
such that either
every $s\in\suc_{q_{j+1}}(t)\sq q_{j+1}$ has the tag
``one-to-one'' or every such an $s$ has a tag ``constant $x_s$.''
In the first case put $V_t=U_t$ and tag $t$
as ``one-to-one.''
In the second case we can find a subset $V_t$ of $U_t$
of size at least
$|U_t|^{1/2}\geq|\suc_{q_{j+1}}(t)|^{1/4}$ such that
the mapping $V_t\ni s\longmapsto x_s\in 2^m$
is either one-to-one or constant equal to $x_t$.
We tag $t$ accordingly and define
\[
q_j=\bigcup\{(q_{j+1})^s\colon
(\exists t\in q_{j+1}\cap\omega^j)(s\in V_t)\}.
\]
This finishes the ``trimming'' construction.

Note that by the construction for every $k\leq j\leq i$ and
$t\in q_j\cap\omega^j$:
\begin{itemize}
\item $q_j\cap\omega^j=q_{j+1}\cap\omega^j$;

\item $\nor_{q_j}(s)=\nor_{q_{j+1}}(s)$ for every
      $s\in q_j\setminus\omega^j$;

\item $\nor_{q_j}(t)=\log_{b_j}\log_{b_j}|V_t|\geq
      \log_{b_j}\log_{b_j}|\suc_{q_{j+1}}(t)|^{\frac{1}{4}}
      =\nor_{q_{j+1}}(t)
      +\log_{b_j}\frac{1}{4}\geq\nor_p(t)-1$;

\item if $t$ has a tag ``constant $x_t$'' then every
      $s\in (q_j)^t\cap\left(\bigcup_{j\leq l\leq i}\omega^l\right)$
      has also the tag ``constant $x_t$.''
\end{itemize}
In particular
$\norsup_{q_k}(u)\geq\norsup_p(u)-1\geq\norsup_p(k)-1\geq n$.
Thus if we put $q_{m,u}=q_k$ then
$\norsup_{q_{m,u}}(u)\geq n$
and either $u$ has a tag ``one-to-one'' or ``constant $x_u$.''
Moreover
in the second case all $s\in q_{m,u}\cap\omega^i$ have the same tag
``constant $x_u$.''

Now if for some $m<\omega$ every $u\in p\cap\omega^k$ is tagged in
$q_{m,u}$ as ``one-to-one'' then
it is easy to see that
\[
q=\bigcup\{q_{m,u}\colon u\in p\cap\omega^k\}
\]
has a representation as in (\ref{eqZZ}).
Indeed, for every $s\in q\cap\omega^i$ let $j_s$ be the
largest $j\leq i$ such that $s\restr j$ is tagged ``one-to-one'' in
$q_{m,u}$. Let $A=\{s\restr j_s\colon s\in q\cap\omega^i\}$.
Then $\bigcup_{s\in A}q^s$ is the required representation.

Thus it is enough to prove that there exist an arbitrarily
large $m$ such that all $u\in p\cap\omega^k$ have a
tag ``one-to-one'' in $q_{m,u}$.

By way of contradiction assume that this is not the case. Then
there exist an infinite set
$X_0\sq\omega$ and $u \in p\cap\omega^k$
such that for every $m\in X_0$ there exists $x_m\in 2^m$
with $u$ having a tag ``constant $x_m$'' in $q_{m,u}$.
In particular,
\[
q_{m,u}\forces \tau\restr m=x_m.
\]
By induction choose an infinite sequence
$X_0\supset X_1\supset X_2\supset\cdots$
of infinite sets such that for every $i<\omega$
there exist $y_i\in 2^i$ and
$T_i\sq\omega^{\leq i}$ with the property that
$q_{m,u}\cap\omega^{\leq i}=T_i$
and $x_m\restr i=y_i$ for every $m\in X_i$.

Choose an infinite set $X=\{m_i<\omega\colon i<\omega\}$
such that $m_i\in X_i$ for every $i<\omega$ and let
$q^\prime=\lim_{i\to\infty} q_{m_i,u}=\bigcup_{i<\omega} T_i$.
Then for every $t\in q^\prime\cap T_i$
we have $\nor_{q^\prime}(t)=\nor_{q_{m_i,u}}(t)\geq \nor_{p}(t)-1$.
Thus $q\in\poset$ and $q\geq p$, as $q_{m_i,u}\geq p$
for every $i<\omega$.
So it is enough to prove that
\begin{equation}\label{eqNM}
\text{
$q^\prime\forces``\tau\restr j=y_j$'' \ \ \ \ \ for every $j<\omega$}
\end{equation}
since then
$y=\bigcup_{j<\omega} y_j\in 2^\omega\cap V$ and
$q^\prime\forces``\tau=y\in V$,'' contradicting the fact that
$p\forces``\tau\notin V$.''

To see (\ref{eqNM}) fix a $j<\omega$ and let $l<\omega$ be such that
$l\geq j$ and $l>i_j$. Take an $m\in X_l\sq X_j$ such that $m\geq j$.
Then $q_{m,u}\cap \omega^{\leq l}=T_l=q^\prime\cap \omega^{\leq l}$.

Fix an arbitrary
$s\in q_{m,u}\cap \omega^l=q^\prime\cap \omega^l$. Then
$q_{m,u}^s\forces``\tau\restr m=x_m$'' while $x_m\restr j=y_j$,
since $m\in X_j$. Thus
\[
q_{m,u}^s\forces\tau\restr j=y_j.
\]
But $s\in q_{m,u}\cap \omega^l\sq p\cap\omega^{>i_j}$.
So there exists an $x_s\in 2^j$ with the property  that
$p^s\forces``\tau\restr j=x_s$.'' Since $q_{m,u}^s\geq p^s$
we conclude that
$q_{m,u}^s$ forces the same thing and so $x_s=y_j$.
Thus,
$p^s\forces``\tau\restr j=y_j$.''
But $(q^\prime)^s\geq p^s$. So
\[
(q^\prime)^s\forces\tau\restr j=y_j
\]
as well. Since it happens for every
$s\in q^\prime\cap \omega^l$ and $j<\omega$ was arbitrary,
we conclude~(\ref{eqNM}).
\qed

The next corollary is equivalent of
Lemma~\ref{lem:Split} for $\alpha=1$. It will
not be used in a sequel. However the same
approach will be used
in the proof of Lemma~\ref{lem:Split} in its general form,
and the proof presented here can shed some light
on what follows.

\cor{lem:P0cor}{Let $p\in\poset$ and $\tau$ be a $\poset$-name
such that
\[
p\forces\tau\in 2^\omega\setminus V
\]
and $p$ reads $\tau$ continuously.
Then for every $n<\omega$
there exist an $m<\omega$, nonempty disjoint sets
$B_0,B_1\subset 2^m$, and $p_0,p_1\geq_n p$
such that
\[
p_i\forces \tau\restr m\in B_i
\]
for $i<2$.
}

Proof. Let $k<\omega$ be such that $\norsup_p(k)\geq n+5$.
Then, by Lemma~\ref{lem:P0},
there exist $m,i_m<\omega$, and $q\geq_{n+4} p$
such that
\[
q=\bigcup_{t\in A}p_t,
\]
where
$A\sq p\cap\left(\omega^{\leq i_m}\setminus\omega^{<k}\right)$,
the elements of $A$ are pairwise incompatible,
$p_t\geq p^t$ for every $t\in A$, and for every $t\in A$
there exists a one-to-one mapping
$h_t\colon \suc_q(t)\to 2^m$ such that
\[
q^s\forces\tau\restr m=h_t(s)
\]
for every $s\in\suc_q(t)$.

Let $\{t_j\colon j<M\}$ be a one-to-one enumeration of $A$
such that $|t_j|\leq|t_{j+1}|$ for every $j<M-1$.
By induction on $j<M$ we will choose a sequence
$\la C^i_j\colon i<2\ \&\ j<M\ra$ such that for every
$i<2$ and $j<M$
\begin{itemize}
\item $C^i_j\in[\suc_q(t_j)]^{(b_l)^{(b_l)^n}}$, where $l=|t_j|$, and

\item the sets $\{h_{t_j}[C^i_j]\subset 2^m\colon i<2\ \&\ j<M\}$
are pairwise disjoint.
\end{itemize}

Given $\la C^i_r\colon i<2\ \&\ r<j\ra$
the choice of $C^0_j$ and $C^1_j$ is possible since for $l=|t_j|$
\[
\left|\bigcup_{i<2,\, r<j}C^i_r\right|\leq 2j\ (b_l)^{(b_l)^n}
\leq 2 \left|T^\star\cap\omega^{\leq l}\right|\ (b_l)^{(b_l)^n}\leq
(b_l)^{(b_l)^{n+2}}
\]
and $|\suc_q(t_j)|\geq (b_l)^{(b_l)^{n+4}}$
so we can choose disjoint $C^0_j,C^1_j\in[\suc_q(t_j)]^{(b_l)^{(b_l)^n}}$
with
\[
h_{t_j}[C^0_j\cup C^1_j]\cap
\left(\bigcup_{i<2,\, r<j}h_{t_r}[C^i_r]\right)=\emptyset.
\]

For $i<2$ define $p_i=\bigcup\{q^s\colon s\in\bigcup_{j<M}C^i_j\}$
and $B_i=\bigcup_{j<M}h_{t_j}[C^i_j]$.
It is easy to see that they have the required properties. \qed

Let us also note the following easy fact.

\lem{lem:P0easy}{Let $\rational$ be an arbitrary forcing, $q\in\rational$,
and let $\tau$ be a $\rational$-name
such that
\[
q\forces\tau\in 2^\omega\setminus V.
\]
Then for every $N<\omega$ there exists
an $m_0<\omega$ with the following property.
If $m_0\leq m<\omega$ then there exist $\{q_n\geq q\colon n<N\}$,
and a
one-to-one sequence
$\la z_n\in 2^m\colon n<N\ra$
such that
\[
q_n\forces \tau\restr m=z_n
\]
for every $n<N$.
}

Proof. By induction on $n<N$ define infinite sequences
$\{x_i^n\in 2^i\colon i<\omega\}$ and
$q\leq q_0^n\leq q_1^n\leq q_2^n,\ldots$
such that for every $i<\omega$
\[
q^n_i\forces \tau\restr i=x^n_i.
\]
Moreover if $x^n=\bigcup_{i<\omega}x^n_i\in 2^\omega\cap V$,
then the construction will be done making sure that
$x^n\notin\{x^k\colon k<n\}$.
It is possible, since
$\{x^k\colon k<n\}\in V$, while
$q$ forces that $\tau$ is not in $V$.

Now choose $m_0<\omega$ such that all restrictions
$\{x^n\restr m_0\colon n<N\}$ are different.
Then for $m_0\leq m<\omega$
define $z_n=x^n\restr m$ and $q_n=p^n_m$ for every $n<N$.
Clearly they have the desired properties. \qed

\rem{rem}{In the text that follows
(including the next lemma) we will often
identify forcing
$\poset_\alpha$ with $\poset_\beta\star\poset_\gamma^\star$,
where $\beta+\gamma=\alpha$ and $\poset_\gamma^\star$
is a $\poset_\beta$-name for $\poset_\gamma$,
via mapping
$\poset_\alpha\ni p\mapsto
\la p\restr\beta,p\restr\alpha\setminus\beta\ra
\in\poset_\beta\star\poset_\gamma^\star$.
However, although this mapping
is an order embedding onto a dense subset of
$\poset_\beta\star\poset_\gamma^\star$, it is not onto.
Thus, each time we will be identifying
an element
$\la p^\prime,q^\prime\ra\in\poset_\beta\star\poset_\gamma^\star$
with a $q\in \poset_\alpha$, in reality
we will be defining $q$ as such an element of $\poset_\alpha$
such that $q\restr\beta\geq_{F,n} p^\prime$
and
$q\restr\beta\forces``q\restr\alpha\setminus\beta=q^\prime$''
for the current values of $F$ and $n$.
To define such a $q$ first
find $q\restr\beta\in\poset_\beta$ and a countable
set $A\sq\alpha$ such
that $q\restr\beta\geq_{F,n} p^\prime$
and $q\restr\beta$ forces that the domain of $q^\prime$
is a subset of $A$. (See~\cite[Lemma~1.6, p.~81]{ShPF}.
Compare also \cite[Lemma~2.3(iii)]{BL}.)
Then it is enough to extend $q\restr\beta$ to
$q\in\poset_\alpha$ with the domain equal to
$A\cup\dom(q\restr\beta)$ in such a way that
$q\restr\xi\forces``q(\xi)=q^\prime(\xi)$'' for every $\xi\in A$.
}

Using Lemma~\ref{lem:P0easy} we can obtain
the following modification of
Lemma~\ref{lem:P0}.
In its statement we will use the symbol
$p|s$ associated with $p\in\poset_\delta$
and $s\in p(0)$ to denote an element of $\poset_\delta$
such that $\dom(p|s)=\dom(p)$, $(p|s)(0)=[p(0)]^s$, and
$(p|s)\restr(\delta\setminus\{0\})=p\restr(\delta\setminus\{0\})$.

\lem{lem:P0iter}{Let $1<\delta<\omega_2$,
$p\in\poset_\delta$, and $\tau$ be a $\poset_\delta$-name
such that
\[
p\forces\tau\in 2^\omega\setminus V[G_1].
\]
Then for every $n,k<\omega$ with $\norsup_{p(0)}(k)\geq n$
there exist an arbitrarily large number $m<\omega$,
$q\geq_{\{0\},n} p$, and a $\poset_1$-name
$\varphi$ such that for every $t\in q(0)\cap\omega^k$
\[
q|t\forces \varphi
\text{ is a one-to-one function from
$\suc_{q(0)}(t)$ into $2^m$}
\]
and
\[
q|s\forces\tau\restr m=\varphi(s)
\]
for every $s\in\suc_{q(0)}(t)$.
}

Proof. Identify $\poset_\delta$ with $\poset_1\star\rational$
and $p$ with $\la p(0),\bar p\ra$, where $\rational$ is a
$\poset_1$-name for $\poset_\gamma$ and $1+\gamma=\delta$.
Let $S=T^\star\cap\omega^{k+1}$ and $N=|S|$.

Take a $V$-generic filter $H$ in $\poset_1$ such that $p(0)\in H$.
For a moment we will work in the model $V[H]$.
In this model let $\tilde\rational$ and $\tilde p$
be the $H$-interpretations of $\rational$ and $\bar p$, respectively.
Moreover let $\tilde\tau\in V[H]$ be a $\tilde\rational$-name
such that $\tilde p$ forces that $\tilde\tau=\tau$.
Then $\tilde p$ forces that $\tilde\tau\in 2^\omega\setminus V[H]$.
Thus,  by Lemma~\ref{lem:P0easy} used in $V[H]$ to $\tilde\tau$,
there exists an  $m_0<\omega$ such that for every $m\geq m_0$
there are
$\{q_s\geq \tilde p\colon s\in S\}$,
and a one-to-one function $f\colon S\to 2^m$
such that
\[
q_s\forces\tilde\tau\restr m=\tau\restr m=f(s)
\]
for every $s\in S$.

Let $\mu$ be a $\poset_1$-name for $m_0$.
Then, by Lemma~\ref{cor00}, there exists $p^\prime\in\poset_1$
and an arbitrarily large $m<\omega$
such that $p^\prime\geq_n p(0)$ and $p^\prime\forces``\mu\leq m$.''

Now let $\{q^\star_s\geq q\colon s\in S\}$
and $\varphi$ be the $\poset_1$-names for
$\{q_s\geq q\colon s\in S\}$
and $f\colon S\to 2^m$, respectively, such that
$p^\prime$ forces the above properties about them.
Moreover let $q^\prime$ be a $\poset_1$-name for
an element of $\rational$ such that
\[
[p(0)]^s\forces q^\prime=q^\star_s
\]
for every $s\in p^\prime\cap \omega^{k+1}$.
Put $q=\la p^\prime,q^\prime\ra$.
It is easy to see that $m$, $q$ and $\varphi$ have the desired
properties. \qed

Lemmas~\ref{lem:P0} and~\ref{lem:P0iter} can be combined together in the
following corollary. Its form is a bit awkward,
but it will allow us to combine two separate cases into
one case in the proof of
Lemma~\ref{lem:Split}.

\cor{cor:P0gen}{Let $1\leq\delta<\omega_2$,
$p\in\poset_\delta$, and $\tau$ be a $\poset_\delta$-name
such that for every $\gamma<\delta$
\[
p\forces\tau\in 2^\omega\setminus V[G_\gamma].
\]
Moreover if $\delta=1$ assume additionally that
$p$ reads $\tau$ continuously.
Then for every $n,k<\omega$ with $\norsup_{p(0)}(k)\geq n+1\geq 2$
there exist an arbitrarily large number $m<\omega$, $i<\omega$
with $i>k$, and $q\geq_{\{0\},n} p$ such that
$q(0)$ can be represented as
\begin{equation}\label{eqZZgen}
q(0)=\bigcup_{t\in A}p_t,
\end{equation}
where
$A\sq p(0)\cap\left(\omega^{\leq i}\setminus\omega^{<k}\right)$ and
the elements of $A$ are pairwise incompatible (as functions).
Moreover
for every $t\in A$ we have $p_t\geq [p(0)]^t$,
and there exists a $\poset_1$-name $\varphi_t$
such that
\[
q|t\forces \varphi_t
\text{ is a one-to-one mapping from $\suc_q(t)$ into $2^m$}
\]
and
\[
q|s\forces\tau\restr m=\varphi(s)
\]
for every $s\in\suc_{q(0)}(t)$.
}

Proof. For $\delta>1$ use
Lemma~\ref{lem:P0iter} with $i=k+1$ and put
$A=q(0)\cap\omega^k$.

For $\delta=1$ use Lemma~\ref{lem:P0}
taking  as $\varphi_t$ the standard names for
the maps
$\suc_q(t)\ni s\longmapsto x_s\in 2^m$.\qed

Next we will consider several properties of the iteration of
forcing $\poset$.

For $p\in\poset_\alpha$, where $\alpha\leq\omega_2$,
and $\sigma\colon F\to\prod_{i<k}n_i\subset\omega^k$,
where $k<\omega$ and $F\in[\alpha]^{<\omega}$,
define a function $p|\sigma$ as follows.
The domain of $p|\sigma$ is equal to $\dom(p)$, and
$(p|\sigma)\restr(\dom(p)\cap\beta)$
is defined by induction on $\beta\leq\alpha$:
\begin{itemize}
\item $(p|\sigma)\restr(\dom(p)\cap\beta)=
      \bigcup_{\gamma<\beta}(p|\sigma)\restr(\dom(p)\cap\gamma)$
      if $\beta$ is a limit ordinal;

\item if $\beta=\gamma+1$ we put
      $(p|\sigma)\restr(\dom(p)\cap\beta)=
      (p|\sigma)\restr(\dom(p)\cap\gamma)$
      provided $\gamma\notin\dom(p)$;

\item if $\beta=\gamma+1$ and $\gamma\in\dom(p)$ we define
      $(p|\sigma)(\gamma)$ as follows:

      \subitem{(A)} if
               $(p|\sigma)\restr(\dom(p)\cap\gamma)\notin\poset_\gamma$
               we define $(p|\sigma)(\gamma)$ arbitrarily;

      \subitem{(B)}
               if $(p|\sigma)\restr(\dom(p)\cap\gamma)\in\poset_\gamma$
               then we put $(p|\sigma)(\gamma)=\tau$ where
               $\tau$ is a $\poset_\gamma$-name such that
               \[
              (p\restr\gamma)|(\sigma\restr\gamma)
              \forces``\tau=[p(\gamma)]^{\sigma(\gamma)}\text{''}
              \]
              if $\gamma\in F$, and $\tau=p(\gamma)$ if $\tau\notin F$.
\end{itemize}
We say that $\sigma$ is {\em consistent}\/ with $p$ if
$p|\sigma$ belongs to $\poset_\alpha$, i.e., when case (A) was never used in
the above definition.
We will be interested in function $p|\sigma$ only
when $\sigma$ is consistent with $p$. In this case intuitively
$p|\sigma$ represents a condition $q\in\poset_\alpha$
with the same domain that $p$ such that
$q(\gamma)=p(\gamma)$ for every $\gamma\not\in F$
and $q(\gamma)=[p(\gamma)]^{\sigma(\gamma)}$ for $\gamma\in F$.
We will use a symbol $\ccc(p,F,k)$ to denote the set of all
$\sigma\colon F\to\omega^k$ consistent with $p$.

Note that if $s\in p(0)$ then function $p|s$ used in
Proposition~\ref{prop2} is equal to
$p|\sigma$, where $\dom(\sigma)=\{0\}$ and $\sigma(0)=s$.
Also such $p|s$ belongs to $\poset_\alpha$ if and only if $s\in p(0)$.
Thus we will identify $\ccc(p,\{0\},k)$ with $p(0)\cap\omega^k$.

For $F\in[\alpha]^{<\omega}$ and $k<\omega$ we say that
$p\in\poset_\alpha$ is {\em $\la F,k\ra$-determined}\/ if for every
$\beta\in F\cap\dom(p)$ and
$\sigma\colon F\cap\beta\to\omega^k$ consistent with $p$
the condition
$(p\restr\beta)|\sigma$
decides already the value of $p(\beta)\cap\omega^k$, that is, if
for every $s\in \omega^k$
\[
\text{either }\ \ \
(p\restr\beta)|\sigma\forces``s\in p(\beta)
\text{''\ \ \   or }\ \ \
(p\restr\beta)|\sigma\forces``s\notin p(\beta).
\text{''}
\]
Note that each $p\in\poset_\alpha$ is $\la\{0\},k\ra$-determined.
%
%
Notice also that for every $p\in\poset_\alpha$,
$k<\omega$, and $F\in[\alpha]^{<\omega}$
if $p$ is $\la F,k\ra$-determined then
\begin{equation}\label{conAnti}
\{p|\sigma\colon\sigma\in\ccc(p,F,k)\}\ \
\text{ is a maximal antichain above }\ p.
\end{equation}
This can be easily proved by induction on $|F|$.
In the same setting we also have
\[
\ccc(p,F\cap\beta,k)=\ccc(p\restr\beta,F\cap\beta,k)=
\{\sigma\restr\beta\colon\sigma\in\ccc(p,F,k)\}
\]
and
\[
(q\restr\beta)|\sigma=(q\restr\beta)|(\sigma\restr\beta)
\]
for every $\beta\leq\alpha$ and $\sigma\in\ccc(p,F,k)$.

\lem{lem10Iter}{Let $\alpha\leq\omega_2$,
$\tau$ be a $\poset_\alpha$-name, $X\in V$ be finite,
and $p\in\poset_\alpha$ be such that
\[
p\forces \tau\in X.
\]
If $i<\omega$ is such that $|X|\leq (b_i)^2$, $t\in p(0)\cap\omega^i$,
and $n<\omega$
is such that $\norsup_{p(0)}(t)\geq n\geq 1$
then there exist $p_t\in\poset_\alpha$ extending $p|t$
and $x\in X$ such that $\norsup_{p_t(0)}(t)\geq n-2$ and
\[
p_t\forces\tau= x.
\]
}

Proof. Let
\[
D=\{T\in\poset\colon(\exists q\geq p|t)
(\exists x\in X)(T=q(0)\ \&\ q\forces``\tau= x\text{''})\}.
\]
Clearly $D$ is dense above $[p(0)]^t$.
We will prove the lemma by induction on~$r^n_D(t)$,
as defined on page~\pageref{page:Rank}.

If $r^n_D(t)=0$ then it is obvious.

If $r^n_D(t)=\alpha>0$ choose $U\in[\suc_{p(0)}(t)]^{\geq (b_i)^{(b_i)^{n-1}}}$
from the definition of $r^n_D(t)$.
By the inductive assumption for every
$s\in U$ there exists $T_s\in D$ extending $[p(0)]^t$ such that
$\norsup_{T_s}(s)\geq n-2$. Choose $q_s$ and $x_s$ witnessing $T_s\in D$,
i.e., such that $q_s\geq p|t$, $q_s(0)=T_s$ and
\[
q_s\forces\tau=x_s.
\]
Since $|X|\leq (b_i)^2$, we can find an
$x\in X$ and $V\sq U$ of cardinality
greater than or equal to $|U|/|X|\geq
(b_i)^{(b_i)^{n-1}}/(b_i)^2=(b_i)^{(b_i)^{n-2}}$
such that $x_s=x$ for every $s\in V$.

Let $S=\bigcup\{q_s(0)\colon s\in V\}$. Then
$\norsup_S(t)=\log_{b_i}|V|\geq n-2$.
Take $p_t\geq p$ such that $\dom(p_t)=\bigcup_{s\in V}\dom(q_s)$,
$p_t(0)=S$, and for $\beta\neq 0$
\[
(p_t\restr\beta)|s\forces``p_t(\beta)=q_s(\beta)\text{''}
\]
for every $s\in V$.
Then $p_t$ satisfies the lemma. \qed

\lem{lem10IterSuper}{Let $\alpha\leq\omega_2$, $p\in\poset_\alpha$,
$k\leq i<\omega$, $\la X_l\colon k\leq l\leq i\ra$ be a sequence
of finite subsets from $V$, and
$\la \tau_l\colon k\leq l\leq i\ra$ a sequence  of \
$\poset_\alpha$-names.
Assume that for every $k\leq l\leq i$
\[
p\forces \tau_l\in X_l
\]
and $|X_l|\leq b_l$. If $n<\omega$ is such that
$\norsup_{p(0)}(k)\geq n+2\geq 3$ then there exist
a family \
$\left\{x_t\in\bigcup_{k\leq l\leq i}X_l\colon t\in p(0)\cap
\left(\bigcup_{k\leq l\leq i}\omega^l\right)\right\}$ \
and $q\geq p$
with the property 
that $p(0)\cap\omega^k=q(0)\cap\omega^k$,
$\norsup_{q(0)}(k)\geq n$,
and
\[
q|t\forces\tau_{|t|}= x_t
\]
for every $t\in p(0)\cap\left(\bigcup_{k\leq l\leq i}\omega^l\right)$.
}

Proof. For every $k\leq l\leq i$ let
$Y_l=\prod_{k\leq j\leq l}X_j$ and notice that
\[
|Y_l|\leq \prod_{k\leq j\leq l}b_j\leq
\left(\prod_{k\leq j< l}n_j\right)\cdot b_l
\leq n_{l-1}!\cdot b_l\leq (b_l)^2.
\]
So,
by Lemma~\ref{lem10Iter}, for every $t\in p(0)\cap\omega^i$
there exist
$p_t\in\poset_\alpha$ extending $p|t$
and $y_t\in Y_i$ such that $\norsup_{p_t(0)}(t)\geq n$ and
\[
p_t\forces\tau_l= y_t(l)
\]
for every $k\leq l\leq i$. We can also assume that all
conditions $p_t$
have the same domain~$D$.

Now let $S_i=p(0)\cap\omega^{\leq i}$. We will construct inductively
a sequence  of trees
$S_i\supset S_{i-1}\supset\cdots\supset S_k$,
such that for every $k\leq l< i$

\pagebreak

\begin{description}
\item{(a)} $S_l\cap\omega^l=S_{l+1}\cap\omega^l$;
\item{(b)}
$\suc_{S_l}(t)=\suc_{S_{l+1}}(t)$ for every $t\in S_l$ with $|t|>l$;
\item{(c)}
$|\suc_{S_l}(t)|\geq (b_l)^{(b_l)^n}$
for every $t\in S_l\cap\omega^l$; and,
\item{(d)} for every $s\in S_l\cap\omega^l$
there exists $y_s\in Y_l$
with the property that
\[
y_t\restr(l+1)=y_s\ \
\text{ for every $t\in S_l\cap\omega^i$ with \ $s\sq t$.}
\]
\end{description}

To make an inductive step take an $l<i$, $i\geq k$, for which $S_{l+1}$
is already defined. For each $s\in S_{l+1}\cap\omega^l$
choose $y_s\in Y_l$ and $L_s\in[\suc_{S_{l+1}}(s)]^{\geq(b_l)^{(b_l)^n}}$
such that
\[
y_s=y_t\restr(l+1)\ \ \text{ for every $t\in L_s$.}
\]
Such a choice can be made, since
$|\suc_{S_{l+1}}(s)|=|\suc_{S_i}(s)|\geq (b_l)^{(b_l)^{n+2}}$
(by the assumption that $\nor_{p(0)}(l)\geq\norsup_{p(0)}(k)\geq n+2$)
while $|Y_l|\leq (b_l)^2$.
Define $L=\bigcup\{L_s\colon s\in S_{l+1}\cap\omega^l\}$
and
\[
S_l=\{s\in S_{l+1}\colon
\text{ either $|s|\leq l$ or $t\sq s$ for some $t\in L$}\}.
\]
This finishes the inductive construction.

Now put $T=\bigcup\{[p(0)]^t\colon t\in S_k\cap\omega^i\}$, and
for every $t\in S_k\cap\left(\bigcup_{k\leq l\leq i}\omega^l\right)$
define $x_t= y_t(|t|)$.
Let $q\in\poset_\alpha$ be such that $\dom(q)=D$,
$q(0)=T$, and $(q\restr\beta)|t\forces``q(\beta)=p_t(\beta)$''
for every $\beta\in D$, $\beta>0$, and $t\in S_k\cap\omega^i$.
It is easy to see that $q$ and all $x_t$'s satisfy the requirements.
\qed

%
%

%

\lem{lemDetGeneral}{Let $\alpha\leq\omega_2$, $k,n<\omega$,
$0\in F\in[\omega_2]^{<\omega}$, and $p\in\poset_\alpha$
be such that
\[
p\restr\beta\forces \norsup_{p(\beta)}(k)\geq n+2\geq 3
\]
for every $\beta\in F$.
Moreover assume that
$k\leq i<\omega$, $\la X_l\colon k\leq l\leq  i\ra$ is a sequence
of finite subsets from $V$, and
$\la \tau_l\colon k\leq l\leq  i\ra$ a sequence  of
$\poset_\alpha$-names with the properties that
for every $k\leq l\leq  i$
\[
p\forces\tau_l\in X_l,
\]
$|X_l|\geq 2$,
and $|X_l|^{(n_{l-1}!)^{2|F|}}\leq b_l$.
Then there exists $q\geq_{F,n}p$ with the following properties.
For every $k\leq l\leq i$
\begin{itemize}
\item $q$ is $\la F,l\ra$-determined; and,
\item there exists a family
$\{x_s\in X_l\colon s\in (\omega^l)^F\ \&\
            s \text{ is consistent with } q\}$
such that
\[
q|s\forces\tau_l= x_s,
\]
for every $s\in(\omega^l)^F$  consistent with $q$.
\end{itemize}
}

Proof.
The proof will be by induction on $m=|F|$.

If $m=|F|=1$ then $F=\{0\}$ and the conclusion follows from
Lemma~\ref{lem10IterSuper}.
(Every $p\in\poset_\alpha$ is $\la\{0\},l\ra$-determined.)

So assume that $m=|F|>1$ and let $\beta=\max F$. Then $0<\beta<\alpha$
and $\poset_\alpha$ is equivalent to
$\poset_\beta\star\poset^\star_\gamma$ where $\beta+\gamma=\alpha$
and $\poset^\star_\gamma$ is a $\poset_\beta$-name
for $\poset_\gamma$.
Let $p_0=\la p\restr\beta,\pi_1\ra\in\poset_\beta\star\poset^\star_\gamma$
be such that $\la p\restr\beta,\pi_1\ra$ is stronger then $p$ and
$p\restr\beta\forces``p(\beta)=\pi_1(0)$.''
Then $p_0\geq_{F,n}p$. Thus we can replace $p$ with~$p_0$.

To make an inductive step, for every $l\leq i$, $l\geq k$, define
\[
X_l^\prime=
\bigcup\left\{(X_l)^T\colon T\sq\prod_{j<l}n_j\subset\omega^l\right\}.
\]
Then
\[
|X_l^\prime|\leq 2^{\left|\prod_{j<l}n_j\right|}\cdot
|X_l|^{\left|\prod_{j<l}n_j\right|}
\leq 2^{n_{l-1}!} |X_l|^{n_{l-1}!}
=(2 |X_l|)^{n_{l-1}!}
\leq |X_l|^{(n_{l-1}!)^2}.
\]
In particular
\[
|X_l^\prime|^{(n_{l-1}!)^{2|F\cap\beta|}}\leq
\left(|X_l|^{(n_{l-1}!)^2}\right)^{(n_{l-1}!)^{2(|F|-1)}}=
|X_l|^{(n_{l-1}!)^{2|F|}}\leq b_l.
\]
So the sequence $\la X_l^\prime\colon k\leq l\leq  i\ra$
and $F\cap\beta$ satisfy the size requirements of the inductive
assumptions.

Now, for a moment, we will work in a model $V[H_\beta]$,
where $H_\beta$ is a $V$-generic filter in $\poset_\beta$
containing $p\restr\beta$.
Let $p_1$ be the valuation of $\pi_1$ in $V[H_\beta]$.
By Lemma~\ref{lem10IterSuper}
there exist
$p^{\prime}\in\poset_\gamma$ extending $p_1$ with
$p^{\prime}(0)\cap\omega^k=p_1(0)\cap\omega^k$
and $\norsup_{p_1(0)}(k)\geq n$,
and for every $l\leq i$, $l\geq k$,
a function
$f_l\colon p^{\prime}(0)\cap\omega^l\to X_l$
such that
\[
p^{\prime}|t\forces\tau_l=f_l(t)
\]
for every $t\in p^{\prime}(0)\cap\omega^l$. Note that,
$f_l\in X_l^\prime$.

Let $\varphi_l$ and $\pi$ be the $\poset_\beta$-names for $f_l$
and $p^{\prime}$, respectively,
such that $p\restr\beta$ forces all the above facts about them.
In particular $p\restr\beta\forces``\varphi_l\in X_l^\prime$''
for all appropriate $l$'s, so,
by the inductive assumption, there exist
$q_0\in\poset_\beta$ and for every $k\leq l\leq i$
a family
\[
\left\{f_s\in X_l\colon s\in \left(\omega^l\right)^{F\cap\beta}\ \&\
            s \text{ is consistent with } p\restr\beta\right\}
\]
such that
$q_0$ is $\la F\cap\beta,l\ra$-determined,
$q_0\geq_{F\cap\beta,n}p\restr\beta$, and
\[
q_0|s\forces\varphi_l= f_s
\]
for every $s\in(\omega^l)^{F\cap\beta}$  consistent with $p\restr\beta$.
In particular every $q_0|s$
decides the value of $\pi(0)\cap\omega^l$, since it is
equal to the domain of $\varphi_l$,
and forces that $\norsup_{\pi(0)}(k)\geq n$.

Let $q=\la q_0,\pi\ra$ and
for every $s\in(\omega^l)^F$  consistent with $p$
define
\[
x_s=f_{s\srestr\beta}(s(\beta)).
\]
It is not difficult to see that it has the required properties. \qed

%

%
%

\bigskip

\noindent
{\bf Proof of Lemma~\ref{lem:Split}.} Let
$\alpha\geq 1$, $p$ and $\tau$ be as in the lemma.

Now for arbitrary $\beta<\alpha$, $\beta\geq 1$, let
$\delta\leq\alpha$ be such that $\beta+\delta=\alpha$.
We will identify $\poset_\alpha$ with
$\poset_\beta\star\poset_\delta^\star$,
where $\poset_\delta^\star$ is a $\poset_\beta$-name for
$\poset_\delta$.
We will also identify
$p$ with $\la p\restr\beta,\pi\ra$.
Upon such identification, we can find a $\poset_\beta$-name
$\tau^\star$ such that
\[
p\restr\beta\forces \tau^\star
\text{ is a name for the same object that $\tau$ is}.
\]
In particular
$p\restr\beta\forces``\pi\forces\tau^\star=\tau$.''

Now if $\alpha$ is a successor ordinal number
put $\alpha=\beta+1$.
In this case $p\restr\beta$ forces that $\pi$ and $\tau^\star$
satisfy the assumptions of the Lemma~\ref{lem:ReadCon2},
so there exists a $\poset_\beta$-name $\pi_0$ such that
\[
p\restr\beta\forces\text{
$\pi_0\geq_n\pi$ and
$\pi$ reads $\tau^\star$ continuously}.
\]
We put $q_0=\la p\restr\beta,\pi_0\ra$
and additionally assume that $\beta\in F$.

If $\alpha$ is a limit ordinal, we put
$\pi_0=\pi$ and $q_0=p$.

Now without loss of generality we can assume that $0\in F$
and $n\geq 1$.
We also put $\beta=\max F$ and fix $q\geq q_0$.

By an easy inductive application of Corollary~\ref{Cor:lem11}
$|F|$-many times we can find $k<\omega$ and
$p^{\prime}\geq_{F,n} q$ such that
\[
p^\prime\restr\gamma\forces
\norsup_{p(\gamma)}(k)\geq n+9
\]
for every $\gamma\in F$. We can also increase $k$, if necessary,
to guarantee that
\begin{equation}\label{eq:kSize}
2|F|+2\leq k. 
\end{equation}
Also since
\[
p^\prime\restr\beta\forces``\pi_0\forces\tau^\star=\tau,\text{''}
\]
$p^\prime\restr\beta$ forces that the assumptions
of Corollary~\ref{cor:P0gen} are satisfied. Thus, applying it
to $\pi_0$, $\tau^\star$, and $k$ defined above,
we can find $\poset_\beta$-names $\mu$, $\rho$, $\pi^{\prime}$,
$\A$, and $\psi$
for $m$, $i$, $q$, $A$ and mapping $A\ni t\mapsto p_t$
respectively,
such that
$p^\prime\restr\beta$ forces
\[
\mu,\rho<\omega\ \ \ \&\ \ \ \pi^{\prime}\geq_{\{0\},n+8}\pi_0\ \ \ \&\ \ \
\pi^{\prime}(0)=\bigcup_{t\in\A}\varphi(t)
\text{ is a representation as in (\ref{eqZZgen}).}
\]
Also, by Corollary~\ref{Cor:lem11}, replacing
$\pi^{\prime}$ with an $\geq_{\{0\}n+8}$-stronger condition, if necessary,
we can assume that there are $m,i<\omega$
such that
\[
p^\prime\restr\beta\forces \mu<m\ \ \ \&\ \ \ \rho<i.
\]
Increasing $i$ and $m$, if necessary, we can also assume that
$m\geq 2$ and
\begin{equation}\label{con:im}
\left|2^m\right|^{(n_{i-1}!)^{2|F|}}< b_i.
\end{equation}

Now notice that we can use Lemma~\ref{lemDetGeneral} to
$p^\prime\restr\beta\in\poset_\beta$,
and the sequences
\[
\la \tau_l\colon k\leq l\leq i\ra=
\left\langle
\pi^\prime(0)\cap\omega^l\colon k\leq l\leq i\right\rangle
\]
and
\[
\la X_l\colon k\leq l\leq i\ra=
\left\langle {\cal P}\left(T^\star\cap\omega^{\leq l}\right)
\colon k\leq l\leq i\right\rangle
\]
since $2|F|+2\leq k$ implies that for every $k\leq l\leq i$
\[
|X_l|^{(n_{l-1}!)^{2|F|}}
\leq \left(2^{(n_{l-1}!)^2}\right)^{(n_{l-1}!)^{2|F|}}
\leq (l+2)^{(n_{l-1}!)^{2|F|+2}}\leq b_l,
\]
where the first inequality is justified by the fact that
$|X_l|\leq 2^{(n_{l-1}!)^2}$, which follows from the following
estimation
\[
\left|T^\star\cap\omega^{\leq l}\right|\leq
\sum_{j<l}n_j!\leq\prod_{j<l}n_j!\leq
n_{l-1}!\prod_{j<l-1}n_j!\leq
n_{l-1}!\prod_{j<l-1}n_{j+1}
\leq (n_{l-1}!)^2.
\]
So
we can find $p^{\prime\prime}\in\poset_\beta$
which is $\la F\cap\beta,l\ra$-determined for each $k\leq l\leq i$,
such that
$p^{\prime\prime}\geq_{F\cap\beta,n+8} p^\prime\restr\beta$, and
that $p^{\prime\prime}|s$ determines the value of
$\pi^\prime(0)\cap\omega^l$ for every
$s\in\left(\omega^l\right)^{F\cap\beta}$ consistent with
$p^{\prime\prime}$.

Next notice also that
$A\cap\omega^{\leq l}\sq T^\star\cap \omega^{\leq l}$.
Thus, the above calculation shows that we can also
use Lemma~\ref{lemDetGeneral} to
$p^{\prime\prime}\in\poset_\beta$,
and the sequences
\[
\la\tau_l\colon k\leq l\leq m\ra=
\la\A\cap\omega^{\leq l}\colon k\leq l\leq i\ra
\]
and
\[
\la X_l\colon k\leq l\leq i\ra=
\la {\cal P}(T^\star\cap\omega^{\leq l})\colon k\leq l\leq i\ra.
\]
So
we can find $p^{\prime\prime\prime}\in\poset_\beta$
such that
$p^{\prime\prime\prime}\geq_{F\cap\beta,n+6} p^{\prime\prime}$, and
that $p^{\prime\prime\prime}|s$ determines the value of
$\A\cap\omega^{\leq l}$ for every
$s\in\left(\omega^l\right)^{F\cap\beta}$ consistent with
$p^{\prime\prime\prime}$.

Now let
$q^1=\la p^{\prime\prime\prime},\pi^{\prime}\ra\in\poset_\alpha$.
Then $q^1\geq_{F,n+6} q$,
\[
q^1\restr\gamma\forces
\norsup_{q^1(\gamma)}(k)\geq n+6
\]
for every $\gamma\in F$, and $q^1$
is $\la F,l\ra$-determined for each $k\leq l\leq i$.
Hence, by the condition (\ref{con:im}),
the assumptions of Lemma~\ref{lemDetGeneral}
are satisfied by $q^1$, and the sequences
$\la\tau_l\colon k\leq l\leq m\ra$
and $\la X_l\colon k\leq l\leq i\ra$, where
$X_i=2^m$,
$\tau_i$ is the restriction to $m$ of the term $\tau$ from
the assumptions of Lemma~\ref{lem:Split},
while for $k\leq l<i$ we put $X_l=2$ and $\tau_l$
a standard name for $0$.
So,
we can find $q^2\geq_{F,n+4} q^1$, which
is still $\la F,l\ra$-determined for each $k\leq l\leq i$,
and a family
$\{x_s\in 2^m\colon s\in (\omega^i)^F\ \&\
            s \text{ is consistent with } q^2\}$
such that
\[
q^2|s\forces\tau\restr m= x_s
\]
for every $s\in(\omega^i)^F$  consistent with $q^2$.
Identify $q^2$ with $\la q^2\restr\beta,\pi^2\ra$ and
note that $q^2\restr\beta$ still forces
that
$\pi^2(0)$ has a
representation as in (\ref{eqZZgen})
and it ``determines'' a big part of this
representation in the sense defined above.
Our final step will be
to ``trim'' $q^2$ (of which we will think as of $\ccc(q^2,F,i)$)
to $q^3$ (identified with $\ccc(q^3,F,i)$)
for which we will be able to repeat
the construction from Corollary~\ref{lem:P0cor}.

For this first note that for every $C\sq\ccc(q^2,F,i)$
there exists a condition $q^2|C$ associated with $q^2$
in a similar way that the condition $q^2|\sigma$ is associated to
$\sigma\in\ccc(q^2,F,i)$.
Also we will consider the elements of $\ccc(q^2,F,i)$
as functions from $i\times F$, where we treat
$i\times F$ as ordered lexicographically by $\leq_{lex}$, and
for $\la l,\gamma\ra\in i\times F$ we define
\[
O(l,\gamma)=\{\la j,\delta\ra\in i\times F\colon
\la j,\delta\ra\leq_{lex}\la l,\gamma\ra\}.
\]

Put $C_0=\ccc(q^2,F,i)$ and let
$\{\la l_j,\gamma_j\ra \colon j\leq r\}$ be a decreasing enumeration of
$(i\setminus k)\times F$ with respect to $\leq_{lex}$.
Note that for every $s\in C_0$ we can associate a tag
``constant $x_s$'' for which
$q^2|s\forces \tau\restr m= x_s$. We will construct by induction
on $j\leq r$ a sequence $C_0\supset C_1\supset\cdots\supset C_r$
such that for every $j\leq r$ and $s\in C_j$ the node
$\la l_j,\gamma_j\ra$ of $s$ is either
tagged ``one-to-one'' (in a sense defined below) or
``constant $x_{s,j}$''
in which case
\[
(q^2|C_j)|s[j]\forces \tau\restr m=x_{s,j},
\]
where $s[j]=s\restr O(l_j,\gamma_j)$.
The above requirement is clearly satisfied for $j=0$, since
$q^2|C_0=q^2$, $s[0]=s$,
and so every $s\in C_0$ is tagged by some constant.
Thus the tag ``one-to-one'' does not appear for $j=0$.
For $j>0$ we will use the tag ``one-to-one'' to $s\in C_j$ if for
$W=\{t\in C_j\colon s[j]\subset t\}$ either
\[
\text{the node
$\la l_{j-1},\gamma_{j-1}\ra$ is
tagged ``one-to-one'' for every $t\in W$}
\]
or for every $t\in W$ the node
$\la l_{j-1},\gamma_{j-1}\ra$ of $t$ is
tagged as a ``constant $x_{t,j-1}$''
and for every $s,t\in W$ if $s[j-1]\neq t[j-1]$ then
$x_{s,j-1}\neq x_{t,j-1}$.
Thus if we think of $C_j$ as of tree $T(C_j)$ being formed
from all $\leq_{lex}$ initial segments of elements of $C_j$,
then the mapping
$\suc_{T(C_j)}(s[j])\ni t[j-1]\mapsto x_{t,j-1}\in 2^m$ is one-to-one.


So assume that for some $0<j\leq r$ the set $C_{j-1}$ is already
constructed. To construct $C_j$ consider first the set
$D=\{s\restr[(l_{j-1}+1)\times(\gamma_{j-1}+1)]\colon s\in C_{j-1}\}$
and note that $D=\ccc(q^2|C_{j-1},F\cap(\gamma_{j-1}+1),l_{j-1}+1)$.
Define
\[
D_0=\{s[j]\restr[(l_{j-1}+1)\times(\gamma_{j-1}+1)]\colon s\in C_{j-1}\}.
\]
Since also
$D_0=\{s\restr\dom(s)\setminus\{\la l_{j-1},\gamma_{j-1}\ra\}\colon s\in D\}$
the elements of $D_0$ are predecessors of those from $D$ in a natural sense.
Now, for every $s_0\in D_0$ let $D_{s_0}$ be the set of all successors
of $s_0$ which belong to $D$, that is, 
$D_{s_0}=\{s\in D\colon s_0\subset s\}$.
In what follows we will describe the method of a choice of
subsets $E_{s_0}$ of $D_{s_0}$.
Then we will define $C_j$ by
\[
C_j=\{s\in C_{j-1}\colon
s\restr[(l_{j-1}+1)\times(\gamma_{j-1}+1)]\in E_{s_0}
\text{ for some }s_0\in D_0\}.
\]
Note that by this definition the norms of $q^2|C_j$ and
$q^2|C_{j-1}$ are the same at every
node of a level $\la l,\gamma\ra$
except for $\la l,\gamma\ra=\la l_{j-1}-1,\gamma_{j-1}\ra\}$,
in which case the norm is controlled by the choice of $E_{s_0}$.

Now to choose sets $E_{s_0}\subset D_{s_0}$ fix an $s_0\in D_0$. We would
like to
look at the tags of elements from $D_{s_0}$
and use the procedure from Corollary~\ref{lem:P0cor} to trim $D_{s_0}$.
However the elements of $D_{s_0}$ do not need to have tags.
Thus we will modify this idea in the following way.
Let $Z_{s_0}=\{s[j]\colon s_0\subset s\in C_{j-1}\}$
and notice that the elements of $Z_{s_0}$ are differed from $s_0$
only by a ``tail'' defined on some pairs $\la l,\gamma\ra$
with $l<l_{j-1}$. Since the possible values of these ``tails''
are already determined by $q^2|s_0$ we have
\[
\left|Z_{s_0}\right|\leq\left|\poset\cap\omega^{l_{j-1}}\right|^{|F|}
\leq\left(n_{l_{j-1}}!\right)^{|F|}.
\]
For $t\in Z_{s_0}$ and $E\subset D_{s_0}$ let 
$E[t]=\{s[j-1]\colon s\restr[(l_{j-1}+1)\times(\gamma_{j-1}+1)]\in E\}$.
Then every element of $E[t]$ has a tag, and we can choose a subset
$E'[t]$ of $E[t]$ of size $\geq |E[t]|^{1/4}$ with either all elements
of $E'[t]$ having the same tag, or all having the tag ``constant'' with
different constant values. Then
$E''[t]=\{s\restr[(l_{j-1}+1)\times(\gamma_{j-1}+1)]\colon s\in E'[t]\}$
is an
$\la E,t\ra$-approximation for $E_{s_0}$.
The actual construction of the set $E_{s_0}$ is obtained by using the
above described operation to all elements $t_1,\ldots,t_p$ of $Z_{s_0}$ one at a
time. More precisely, we  put $E_0=D_{s_0}$ and define $E_\nu$ for $1\leq
\nu\leq p$
as $E_{\nu-1}''[t_\nu]$. Then we put $E_{s_0}=E_p$
and note that
\[
|E_{s_0}|\geq |D_{s_0}|^{4^{-|Z_{s_0}|}}\geq
|D_{s_0}|^{4^{-\left(n_{l_{j-1}}!\right)^{|F|}}}\geq
|D_{s_0}|^{\left(b_{l_{j-1}}\right)^{-1}}.
\]
This finishes the inductive construction.

Now define $q^3=q^2|C_r$ and notice that $q^3\geq_{F,n+3} q^2$.
Indeed, this follows from the norm preservation remark above
and the fact that
\[
|E_{s_0}|\geq |D_{s_0}|^{\left(b_{l_{j-1}}\right)^{-1}}\geq
\left((b_{l_{j-1}})^{(b_{l_{j-1}})^{n+4}}\right)^{\left(b_{l_{j-1}}\right)^{-1}}
=(b_{l_{j-1}})^{(b_{l_{j-1}})^{n+3}}.
\]


By the above construction for every $s\in C_r$ every node
$s[j]$ of $s$ from level $\la l_j,\gamma_j\ra$ has a tag
in $q^3$. Moreover, although $s=s[0]$ has a tag ``constant,''
all this tags cannot be ``constant.''  Indeed, if
$l$ is such that
$(q^2\restr\beta)|(s\restr i\times\beta)$
forces that the node $q^2(\beta)(l)$ is tagged ``one-to-one''
while its successors are tagged as constants, then it is
easy to see that the same node (more precisely, the node from level
$\la \max(F\cap \beta),l+1\ra$) will remain tagged
``one-to-one'' in our recent tagging procedure.
In particular, for every $s\in C_r$ there exists a maximal
number $j_s<r$ for which $s[j_s]$ is marked ``one-to-one.''

To make the final step let $T_1=T(C_r)$ be the tree as defined above and let
$\{t_j\colon j<M\}$ be a one-to-one enumeration of
$\{s[j_s]\colon s\in C_r\}$
such that
$|\suc_{T_1}(t_j)|\leq|\suc_{T_1}(t_{j+1})|$ for every $j<M-1$.
We will proceed as in Corollary~\ref{lem:P0cor}.
By induction on $j<M$ we will choose a sequence
$\la C^u_j\colon u<2\ \&\ j<M\ra$ such that for every
$u<2$ and $j<M$
\begin{itemize}
\item if $|\suc_{T_1}(t_j)|=(b_l)^{(b_l)^{n+3}}$ then
$C^i_j\in[\suc_q(t_j)]^{(b_l)^{(b_l)^n}}$; and,

\item the sets $\{h[C^u_j]\subset 2^m\colon u<2\ \&\ j<M\}$
are pairwise disjoint.
\end{itemize}

Given $\la C^u_r\colon u<2\ \&\ r<j\ra$
we can choose  $C^0_j$ and $C^1_j$ since for $l=|t_j|$
\[
\left|\bigcup_{u<2,\, r<j}C^u_r\right|\leq
2j\ (b_l)^{(b_l)^n}
\leq 2 \left|T^\star\cap\omega^{\leq l}\right|^{|F|}\ (b_l)^{(b_l)^n}
\leq (b_l)^{(b_l)^{n+2}}
\]
and $|\suc_{T_1}(t_j)|=(b_l)^{(b_l)^{n+3}}$
therefore
it is possible to choose disjoint sets
$C^0_j,C^1_j\in[\suc_{T_1}(t_j)]^{(b_l)^{(b_l)^n}}$
with
\[
h[C^0_j\cup C^1_j]\cap
\left(\bigcup_{u<2,\, r<j}h[C^u_r]\right)=\emptyset.
\]

For $u<2$ define
$C_u=\bigcup\{(T_1)^s\colon s\in\bigcup_{j<M}C^u_j\}$
and $B_u=\bigcup_{j<M}h[C^u_j]$.
It is easy to see that $p_u=q^2|C_u$ and $B_u$
have the required properties. \qed

\end{document}